\documentclass[10pt, leqno]{article}
\usepackage{amsmath,amsfonts,amsthm,amscd,amssymb}
\numberwithin{equation}{section}
\theoremstyle{plain}
\newtheorem{theo}{Theorem}[section]
\newtheorem{prop}[theo]{Proposition}
\newtheorem{coro}[theo]{Corollary}
\newtheorem{lemm}[theo]{Lemma}
\theoremstyle{definition}
\newtheorem{defi}[theo]{Definition}
\newtheorem{rema}[theo]{Remark}
\newtheorem{theo-defi}[theo]{Theorem-Definition}
\newtheorem{prop-defi}[theo]{Proposition-Definition}
\newtheorem{rem-defi}[theo]{Remark-Definition}

\newtheorem{prob}[theo]{Problem}

\def \be{\beta}
\def \bul{\bullet}
\def \col{\colon}

\def \del{\delta}
\def \eps{\epsilon}
\def \Gam{\Gamma}
\def \gam{\gamma}

\def \kap{\kappa}
\def \Lam{\Lambda}
\def \lam{\lambda}

\def \Lo{\Longrightarrow}
\def \lo{\longrightarrow}
\def \lom{\longmapsto}
\def \mab{\mathbb}

\def \Om{\Omega}

\def \ol{\overline}
\def \os{\overset}
\def \parno{\par\noindent}
\def \part{\partial}

\def \sus{\subset}

\def \ul{\underline}
\def \us{\underset}

\def \vil{\varinjlim}

\def \wh{\widehat}
\def \wt{\widetilde}
\bigskip

\newcommand{\getsfrom}{\ensuremath{
\longleftarrow\kern-.52em\lower-.1ex\hbox%
{$\shortmid\,$}}}

\begin{document}
\title{Kato-Nakayama's comparison theorem and 
analytic log etale topoi}
\author{Yukiyoshi Nakkajima
\date{}\thanks{ \endgraf}}
\maketitle

\parno
\bigskip
\parno 
{\bf Abstract}: In this paper we study topics related to 
one of Kato-Nakayama's comparison theorems in \cite{kn} 
using analytic log etale topoi.  
\bigskip
\parno

\section{Introduction} 
In \cite{kn} 
Kato and Nakayama have proved 
the following comparison theorem (cf.~\cite[(5.9)]{illl}):  
\begin{theo}[{\cite[(0.2) (1); log etale vs.~log Betti]{kn}}]
\label{theo:cole}
Let $X$ be an fs$($=fine and saturated$)$ 
log scheme over ${\mab C}$ 
whose underlying scheme is 
locally of finite type over ${\mab C}$. 
Let $X^{\log}_{\rm an}$ be the real blow up 
of the analytification $X_{\rm an}$ of 
$X$ $(${\rm \cite{kn}}$)$.  
Let 
${\rm D}^+_{{\rm c}{\textrm -}{\rm tor}}(X^{\log}_{\rm et})$ 
be the derived category of bounded below complexes of 
abelian sheaves in $\wt{X}^{\log}_{\rm et}$ 
whose cohomology sheaves are constructible 
{\rm (\cite{kn})} and torsion. 
Let $K^{\bul}$ be an object of 
${\rm D}^+_{{\rm c}{\textrm -}{\rm tor}}(X^{\log}_{\rm et})$.
Let $K^{\bul \log}_{\rm an}$ be the inverse image of 
$K^{\bul}$ in $\wt{X}^{\log}_{\rm an}$. 
Then there exists a canonical isomorphism 
\begin{equation*}
H^h(X^{\log}_{\rm et}, K^{\bul}) \os{\sim}{\lo}
H^h(X^{\log}_{\rm an}, K^{\bul \log}_{\rm an}) 
\quad (h \in {\mab Z}). 
\end{equation*}\label{eqn:ctk}
\end{theo}
\par 
In this paper, for an fs log analytic space $Y$ 
over ${\mab C}$, we introduce a new topos 
$\wt{Y}^{\log}_{\rm et}$ 
of $Y$ ($\wt{Y}^{\log}_{\rm et}$ 
is an analytic analogue of $\wt{X}^{\log}_{\rm et}$) 
and we prove 
the following: 

\begin{theo}\label{theo:logint}
$(1)$ $(${\rm {\bf analytically log etale vs.~log Betti}}$)$ 
Let 
${\rm D}^+_{{\rm lcl}{\textrm -}{\rm tor}}
(Y^{\log}_{\rm et})$ 
be the derived category of bounded below complexes of 
abelian sheaves in $\wt{Y}^{\log}_{\rm et}$ 
whose cohomology sheaves are 
locally classical and torsion 
$($see $(\ref{defi:alelb})$ below for the definition of 
a locally classical abelian sheaf$)$.
Let $K^{\bul}$ be an object of 
${\rm D}^+_{{\rm lcl}{\textrm -}{\rm tor}}
(Y^{\log}_{\rm et})$.  
Let $K^{\bul \log}$ be the inverse image of $K^{\bul}$ 
in $\wt{Y}^{\log}$. 
Then there exists a canonical isomorphism 
\begin{equation*}
H^h(Y^{\log}_{\rm et}, K^{\bul}) \os{\sim}{\lo}
H^h(Y^{\log}, K^{\bul \log}) 
\quad (h \in {\mab Z}). 
\end{equation*}\label{eqn:lectk} 
\par
$(2)$ $(${\rm {\bf GAGA: algebraically log etale 
vs.~analytically log etale}}$)$ 
Let the notations be as in {\rm (\ref{theo:cole})}. 
Let $K^{\bul}_{\rm an}$ be the inverse image of $K^{\bul}$ 
in $\wt{(X_{\rm an})}{}^{\log}_{\rm et}$. 
Then there exists a canonical isomorphism 
\begin{equation*}
H^h(X^{\log}_{\rm et}, K^{\bul}) \os{\sim}{\lo}
H^h((X_{\rm an})^{\log}_{\rm et}, K^{\bul}_{\rm an}) 
\quad (h \in {\mab Z}). 
\end{equation*}\label{eqn:anctk} 
\end{theo} 

\parno
Let 
${\rm D}^+_{{\rm c}{\textrm -}{\rm tor}}(Y^{\log}_{\rm et})$ 
be the derived category of bounded below complexes of 
abelian sheaves in $\wt{Y}^{\log}_{\rm et}$ 
whose cohomology sheaves are constructible and torsion 
$($see {\rm (\ref{defi:ancon})} below for the definition 
of a constructible abelian sheaf in 
$\wt{Y}^{\log}_{\rm et})$. 
Then we shall see that we have a natural functor 
${\rm D}^+_{{\rm c}{\textrm -}{\rm tor}}(X^{\log}_{\rm et})
\lo 
{\rm D}^+_{{\rm c}{\textrm -}{\rm tor}}
((X_{\rm an})^{\log}_{\rm et})$ 
(if $\os{\circ}{X}$ is quasi-compact and quasi-separated) 
and  an inclusion  
${\rm D}^+_{{\rm c}{\textrm -}{\rm tor}}
(Y^{\log}_{\rm et}) \subset 
{\rm D}^+_{{\rm lcl}{\textrm -}{\rm tor}}
(Y^{\log}_{\rm et})$. 
Consequently, 
by (\ref{theo:logint}) (1) and (2),  
we immediately reobtain (\ref{theo:cole}), 
which is another proof of \cite[(0.2) (1)]{kn}.  
\par 
The theorem (\ref{theo:logint}) (1) prompts us to say 
roughly that 
the abstract topos $\wt{Y}^{\log}_{\rm et}$ can replace  
the topological space $Y^{\log}$  
and conversely that the concrete topological space 
$Y^{\log}$ represents the topos 
$\wt{Y}^{\log}_{\rm et}$. The theorem (\ref{theo:logint}) (2) 
is a theorem of GAGA type. 
Though (\ref{theo:logint}) (1) and (2) tempt us to say 
that $\wt{X}^{\log}_{\rm et}$,  
$\wt{(X_{\rm an})}{}^{\log}_{\rm et}$ 
and $\wt{X}^{\log}_{\rm an}$ are the same topoi for 
the calculations of the cohomologies of 
bounded below complexes of abelian sheaves 
whose cohomology sheaves are constructible and torsion, 
I think that $\wt{(X_{\rm an})}{}^{\log}_{\rm et}$ 
is more closely connected with $\wt{X}^{\log}_{\rm et}$ 
than $\wt{X}^{\log}_{\rm an}$ in general 
by taking (\ref{rema:lslt}) (1) below 
(see also (\ref{rema:lslt}) (3) and 
(\ref{theo:lgr}) below)
into account. 
\par 
The contents of this paper are as follows. 
\par 
In \S\ref{sec:lgelk} 
we introduce $\wt{Y}^{\log}_{\rm et}$. 
The topos $\wt{Y}^{\log}_{\rm et}$ plays 
a key role in almost all parts of this paper.  
In particular, a key abelian sheaf 
${\cal M}_{Y,\log}^{\rm gp}$ 
defined in \S\ref{sec:lgelk} 
(${\cal M}_{Y,\log}^{\rm gp}$ is an analytic 
analogue of ${\cal M}_{X,\log}^{\rm gp}$ 
defined in \cite{kn}) lives in $\wt{Y}^{\log}_{\rm et}$ 
(not in $\wt{Y}^{\log}$). 
In \S\ref{sec:knctsc} we give a proof of 
(\ref{theo:logint}) (1) 
by using the analytic log Kummer sequence;  
in \S\ref{sec:bcg} we give a proof of 
(\ref{theo:logint}) (2). 
In fact we prove a base change theorem which is 
a generalization of (\ref{theo:logint}) (2). 
Motivated by (\ref{theo:logint}) (2),  
we give a logarithmic version 
of Grauert-Remmert's theorem 
in the end of \S\ref{sec:bcg}. 
In \S\ref{sec:knct} and \S\ref{sec:les}, 
using ${\cal M}_{X_{\rm an},\log}^{\rm gp}$, 
we obtain a commutative diagram 
(\ref{cd:cogokugcret}) below which 
compares three calculations of 
certain two higher direct images  
by the use of the analytic and algebraic 
log Kummer sequences and 
by the use of the log exponential sequence in \cite{kn}.  
The starting purpose in writing this paper 
was to give the commutative diagram 
because it is necessary in \cite{ndeg} for 
the determination of the delicate sign 
before the \v{C}ech-Gysin morphism appearing in 
the boundary morphism of the $E_1$-terms of  
the $l$-adic weight spectral sequence in \cite{nd} 
from the analogous determination in \cite{ndeg} for 
the $\infty$-adic weight spectral sequence 
essentially obtained in \cite{fn}. 
\bigskip
\parno
{\bf Acknowledgment.}
I am very thankful to C.~Nakayama  
for explaining to me
the detail of an argument in 
the proof of \cite[(2.6)]{kn}. 
Without his explanation, I would misunderstand it.
\bigskip
\par\noindent 
{\bf Notation.} (1) 
For a log scheme 
(resp.~log analytic space over ${\mab C}$) 
$X$ in the sense of Fontaine-Illusie-Kato 
(\cite{klog1}, (resp.~\cite{kn})), 
we denote by 
$\os{\circ}{X}$ 
the underlying scheme (resp.~underlying analytic space) 
of $X$ and by ${\cal M}_X$ the log structure of $X$. 
For a morphism $f \col X \lo Y$ of log schemes 
(resp.~log analytic spaces over ${\mab C}$), 
$\os{\circ}{f} \col \os{\circ}{X} \lo \os{\circ}{Y}$ 
denotes the underlying morphism of $f$ between  
schemes (resp.~analytic spaces over ${\mab C}$).  
For a scheme (resp.~analytic space over ${\mab C}$) 
$Y$, we denote simply by $Y$ 
the trivial log scheme 
(resp.~trivial log analytic space) $(Y,{\cal O}_Y^*)$. 
\par
(2) The word ``etale'' in this paper 
means the word ``\'{e}tale'' in French. 
\bigskip
\parno
{\bf Conventions.}
(1) (\cite[0.3.2]{bbm}) 
Let ${\cal C}$ be an exact additive category.  
For a short exact sequence 
$$0\lo (E^{\bul},d^{\bul}_E) 
\os{f}{\lo} (F^{\bul},d^{\bul}_F) 
\os{g}{\lo} 
(G^{\bul},d^{\bul}_G) \lo 0$$
of bounded below complexes of objects in ${\cal C}$, let 
$(E^{\bul}[1],d^{\bul}_E[1])\oplus 
(F^{\bul},d^{\bul}_F)$ be the mapping cone of $f$.
We fix an isomorphism 
``$(E^{\bul}[1],d^{\bul}_E[1])\oplus 
(F^{\bul},d^{\bul}_F)\owns (x,y) \lom g(y)\in 
(G^{\bul},d^{\bul}_G)$'' 
in the derived category ${\rm D}^+({\cal C})$.
\smallskip
\par
(2) (\cite[0.3.2]{bbm}) 
Under the situation (1), the boundary morphism
$(G^{\bul},d^{\bul}_G) \lo (E^{\bul}[1],d^{\bul}_E[1])$ 
in ${\rm D}^+({\cal C})$ 
is the following composite morphism 
$$(G^{\bul},d^{\bul}_G) 
\os{\sim}{\longleftarrow} (E^{\bul}[1],d^{\bul}_E[1]) \oplus 
(F^{\bul},d^{\bul}_F) 
\os{{\rm proj}.}{\lo} (E^{\bul}[1],d^{\bul}_E[1])
\os{(-1)\times}{\lo} (E^{\bul}[1],d^{\bul}_E[1]).$$
\par 
(3) Let $({\cal T}, {\cal A})$ be a ringed topos. 
Let $({\cal T}_{\bul}, {\cal A}^{\bul})
:=
({\cal T}_{t}, {\cal A}^{t})_{t\in {\mab N}}$ 
be a constant simplicial ringed topos 
defined by $({\cal T}, {\cal A})$: 
${\cal T}_{t}={\cal T}$,
${\cal A}^{t}={\cal A}$.
Let $M$ be a complex of ${\cal A}^{\bul}$-modules.  
The complex $M$ defines a double complex 
$M^{\bul \bul}=(M^{t \bul})_{t\in {\mab N}}$ 
of ${\cal A}$-modules whose boundary 
morphisms will be fixed 
in (\ref{eqn:bdsignss}) below 
(Our convention on the place of cosimplicial degrees 
is different from that in \cite[(5.1.9) IV]{dh3}.).
Let ${\bf s}(M)$ be the single complex 
$\bigoplus_{t+s=n}M^{ts}$
with the following boundary morphism:
\begin{equation*} 
d(x^{ts}) =
\sum_{i= 0}^{t+1}(-1)^{i}
\del_{i}(x^{ts})+
(-1)^{t}d_M(x^{ts}) \quad (x^{ts}\in M^{ts}),
\tag{1.2.1}\label{eqn:bdsignss}
\end{equation*}
where 
$d_M \col M^{ts} \lo  M^{t,s+1}$ 
is the boundary morphism arising 
from the boundary 
morphism of the complex $M$ and 
$\del_{i} \col M^{t s} 
\lo M^{t+1, s}$ 
$(0 \leq i \leq t+1)$ is a 
standard coface morphism. 
The convention on  signs in (\ref{eqn:bdsignss}) is 
different from that in 
\cite[(5.1.9.2)]{dh3}. 

\section{Analytic log etale topoi}\label{sec:lgelk} 
First we recall a well-known method 
(\cite[XI 4]{sga4-3}) quickly to fix our ideas. 
\par
Let $T$ be a topological space. 
Let $T_{\rm cl}$ be a site 
defined by the following: 
\par
(2.0.1) An object of $T_{\rm cl}$ is a 
local isomorphism $U \lo T$ of topological spaces. 
\par
(2.0.2) A morphism in $T_{\rm cl}$ is 
a morphism of topological spaces over $T$.
\par
(2.0.3) A family $\{U_{\lam} \lo U\}_{\lam}$ 
of morphisms in $T_{\rm cl}$ 
is called a covering 
if the union of the images of 
$U_{\lam}$'s is $U$; 
the coverings define 
a Grothendieck pretopology 
and hence a Grothendieck topology on the category 
$T_{\rm cl}$.  
\par
Let $Y$ be an fs log analytic space over ${\mab C}$ 
in the sense of \cite[\S1]{kn}. 
Let $Y^{\log}$ be the real blow up of $Y$ 
(\cite[(1.2)]{kn}). 
Let $Y^{\log}_{\rm cl}$ be the site above 
for the topological space $Y^{\log}$. 
Let $\wt{Y}{}^{\log}$ and $\wt{\os{\circ}{Y}}$ be 
the topoi defined by the classical topologies 
of $Y^{\log}$ and $\os{\circ}{Y}$, respectively. 
Then we have natural morphisms   
$\eps_{\rm cl} \col 
\wt{Y}^{\log}_{\rm cl} \lo \wt{\os{\circ}{Y}}_{\rm cl}$, 
$\eps_{\rm top} \col 
\wt{Y}^{\log} \lo \wt{\os{\circ}{Y}}$, 
$\mu^{\log} \col \wt{Y}^{\log}_{\rm cl} 
\lo \wt{Y}^{\log}$ and  
$\mu \col \wt{\os{\circ}{Y}}_{\rm cl} \lo 
\wt{\os{\circ}{Y}}$  
of topoi fitting into 
the following commutative diagram 
\begin{equation*}
\begin{CD}
\wt{Y}^{\log}_{\rm cl} @>{\mu^{\log}}>> 
\wt{Y}^{\log} \\
@V{\eps_{\rm cl}}VV @VV{\eps_{\rm top}}V \\
\wt{\os{\circ}{Y}}_{\rm cl} @>{\mu}>> \wt{\os{\circ}{Y}}. 
\end{CD} 
\tag{2.0.4}\label{cd:natez}
\end{equation*} 
We sometimes denote $\eps_{\rm cl}$ by $\eps_{Y^{\log}}$. 
\par 
Let ${\cal M}_{Y_{\rm cl}}:=\mu^{-1}({\cal M}_Y)$ 
(resp.~${\cal O}_{Y_{\rm cl}}:=
\mu^{-1}({\cal O}_Y)$) 
be the sheaf of log structures (resp.~the structure sheaf)
in $\wt{\os{\circ}{Y}}_{\rm cl}$.  
Henceforth, in this paper, 
we consider $\wt{\os{\circ}{Y}}_{\rm cl}$, 
$\wt{Y}^{\log}_{\rm cl}$, ${\cal M}_{Y_{\rm cl}}$ 
and ${\cal O}_{Y_{\rm cl}}$ in almost all cases, 
and we denote them simply by 
$\wt{\os{\circ}{Y}}$, 
$\wt{Y}^{\log}$, ${\cal M}_{Y}$ 
and ${\cal O}_{Y}$, respectively, 
as in \cite{kn} unless stated otherwise. 
\par 
Next, let us introduce a topos $\wt{Y}^{\log}_{\rm et}$. 
\par 
As in \cite[(3.1)]{klog1}, 
we can define the exact closed immersion 
of fine log analytic spaces over ${\mab C}$.  
Let 
\begin{equation*} 
\begin{CD}
T_0 @>>> U \\ 
@V{\cap}VV @VV{f}V \\ 
T@ >>> V 
\end{CD} 
\tag{2.0.5}\label{cd:t0tuv}
\end{equation*}
be a commutative diagram of 
fine log analytic spaces over ${\mab C}$, 
where the left vertical morphism is 
an exact closed immersion defined by 
a nilpotent quasi-coherent ideal sheaf of ${\cal O}_T$.  
As in \cite[(3.3)]{klog1}, we say that $f$ is 
{\it log smooth} 
(resp.~{\it log etale}) 
if there exists a (resp.~unique) morphism 
$T\lo U$ locally 
which makes the resulting two diagrams commutative. 
(In the (log) analytic case, 
we do not assume that 
$\os{\circ}{f}$ is locally of finite presentation.) 
As usual, the properties of 
the log smoothness and the log etaleness 
are stable under the composition 
of morphisms and the base change 
in the category of fine log analytic spaces over ${\mab C}$ 
(We can prove that the fiber product exists in  
the category of fine log analytic spaces over ${\mab C}$
using the classical result 
of the existence of the fiber product 
in the category of analytic spaces over ${\mab C}$; 
we can also prove the existence of the fiber product 
in the category of fs log analytic spaces 
over ${\mab C}$.). 
\par
The proof below for local descriptions 
of log smooth and etale morphisms 
in the (log) analytic case are 
slightly different from 
the proof of \cite[(3.5)]{klog1} 
in the log algebraic case 
(we need to give care to the convergence): 

\begin{prop}\label{prop:locdse} 
Let $f \col U \lo V$ be a morphism of 
fine log analytic spaces over ${\mab C}$. 
Then the following conditions $(1)$ and $(2)$ 
are equivalent$:$
\par 
$(1)$ $f$ is log smooth $($resp.~log etale$)$. 
\par 
$(2)$ There exists a local chart 
$Q \lo P$ of $f$ 
satisfying the following conditions 
$(a)$ and $(b):$
\par 
$(a)$ the morphism $Q^{\rm gp}\otimes_{\mab Z}{\mab Q} 
\lo P^{\rm gp}\otimes_{\mab Z}{\mab Q}$ is injective 
$($resp.~isomorphic$)$. 
\par 
$(b)$ The locally induced morphism $\os{\circ}{U} \lo 
\os{\circ}{V}
\times_{{\rm Spec}({\mab C}[Q])_{\rm an}}
{\rm Spec}({\mab C}[P])_{\rm an}$ 
is locally isomorphic as analytic spaces over ${\mab C}$. 
\end{prop}
\begin{proof}
(2)$\Lo$(1): 
Because the morphism $U \lo V
\times_{({\rm Spec}({\mab C}[Q])_{\rm an},Q^a)}
({\rm Spec}({\mab C}[P])_{\rm an},P^a)$ 
is log etale in our sense 
by the easy implication of (\ref{lemm:tse}) (2) below,  
we have only to prove that 
the morphism 
$({\rm Spec}({\mab C}[P])_{\rm an},P^a)\lo
({\rm Spec}({\mab C}[Q])_{\rm an},Q^a)$ is 
log smooth (resp.~log etale).  
Consider the commutative diagram 
(\ref{cd:t0tuv}) for this morphism and 
let $t_0 \col T_0 \lo ({\rm Spec}({\mab C}[P])_{\rm an},P^a)$ 
be the morphism. 
By the same  proof as 
that of \cite[(3.4)]{klog1}, 
we have an (resp.~unique) extension homomorphism 
$P\lo {\cal M}_T$ 
of $P\lo {\cal M}_{T_0}$ 
and  $Q\lo {\cal M}_T$. 
For a point $x \in \os{\circ}{T}_0$,  
let ${\mathfrak m}_x$ 
and ${\mathfrak m}_{0,x}$ 
be the maximal ideals of 
${\cal O}_{T,x}$ and ${\cal O}_{T_0,x}$, respectively. 
Since we have the composite morphism 
$\os{\circ}{T}_0 \os{t_0}{\lo} {\rm Spec}({\mab C}[P])_{\rm an}
\os{\pi}{\lo} {\rm Spec}({\mab C}[P])$ of 
local ringed spaces 
and since ${\cal O}_{T,x}/{\mathfrak m}_x=
{\cal O}_{T_0,x}/{\mathfrak m}_{0,x}$,  
the induced morphism 
$(\pi t_0)^{-1}({\mab C}[P]) \lo {\cal O}_T$ 
by the morphism $P \lo {\cal M}_T$ 
extends to a morphism 
$\os{\circ}{T} \lo{\rm Spec}({\mab C}[P])_{\rm an}$ 
by the universality of the analytification 
(\cite[XII (1.1)]{sga1}). 
Consequently 
we have a (resp.~unique) desired morphism 
$T\lo ({\rm Spec}({\mab C}[P])_{\rm an},P^a)$. 
\par 
(1)$\Lo$(2): Assume that $f$ is log smooth. 
We may assume that there exists a global chart 
$Q \lo {\cal O}_V$ of $V$. 
\par 
Let $\Lam^1_{U/V}$ be the sheaf of 
log differential forms on $U/V$ 
defined similarly in \cite[(1.7)]{klog1} 
(see also \cite{fr1} for the semistable case and 
\cite[(3.5)]{kn} for the absolute case). 
Let ${\cal L}$ be an ${\cal O}_U$-module of finite type. 
Let $\vert U \vert $ be 
the underlying topological space of $U$. 
Let $\os{\circ}{U}({\cal L})
:=(\vert U\vert , {\cal O}_U\oplus {\cal L})$ 
be an analytic space over ${\mab C}$, where we endow 
${\cal O}_U\oplus {\cal L}$ with a natural 
structure of a sheaf of rings by defining 
${\cal L}^2=0$. 
Let $\iota \col \os{\circ}{U} \os{\sus}{\lo} 
\os{\circ}{U}({\cal L})$ be a closed immersion 
and endow $\os{\circ}{U}({\cal L})$ 
with the log structure $\iota_*({\cal M}_U)$. 
Let $U({\cal L})$ be the resulting 
fine log analytic space over 
${\mab C}$.
For a surjective morphism 
${\cal L} \lo {\cal L}_0$ of ${\cal O}_U$-modules 
of finite type, 
we have a closed immersion 
$U({\cal L}_0)\os{\sus}{\lo} U({\cal L})$ 
of fine log analytic spaces over ${\mab C}$.  
Then, using the definition of the smoothness 
and the standard deformation theory 
(cf.~\cite[III (5.1)]{sga1}), 
we easily see that the natural morphism 
${\rm Hom}_{{\cal O}_U}(\Lam^1_{U/V}, {\cal L}) 
\lo {\rm Hom}_{{\cal O}_U}(\Lam^1_{U/V}, {\cal L}_0)$ 
is locally surjective. 
Hence $\Lam^1_{U/V}$ is a locally  
projective ${\cal O}_U$-module of finite type. 
Let $x \in \os{\circ}{U}$ be a point and 
let $d\log t_1,\ldots, d\log t_r$ 
$(t_1, \ldots, t_r\in {\cal M}_{U,x}, r\in {\mab N})$ 
be a basis of $\Lam^1_{U/V,x}$ (If $f$ is log etale, 
we easily see that $\Lam^1_{U/V}=0$ by the standard log 
deformation theory (cf.~\cite[III (5.1)]{sga1}).).  
Then, as in the proof of \cite[(3.5)]{klog1}, 
by using a natural morphism 
${\mab N}^r\oplus Q \lo {\cal M}_{U,x}$ and 
a well-defined surjective morphism 
$$\Lam^1_{U/V,x}\owns d\log a \lom 1\otimes a 
\in \kap(x) \otimes_{\mab Z}
({\cal M}^{\rm gp}_{U,x}/{\cal O}^*_{U,x}
{\rm Im}(\os{\circ}{f}{}^{-1}
({\cal M}_{V,\os{\circ}{f}(x)})))
\quad (a \in {\cal M}_{U,x}),$$
we have a surjective morphism 
${\mab C}\otimes_{\mab Z}({\mab Z}^r\oplus Q^{\rm gp}) 
\lo {\mab C}\otimes_{\mab Z}
({\cal M}^{\rm gp}_{U,x}/{\cal O}^*_{U,x})$. 
Because ${\cal M}^{\rm gp}_{U,x}/{\cal O}^*_{U,x}$ 
is torsion-free, there exists 
a surjective homomorphism $A \lo 
{\cal M}^{\rm gp}_{U,x}/{\cal O}^*_{U,x}$, where 
$A$ is an abelian group 
which has ${\mab Z}^r\oplus Q^{\rm gp}$ 
as a subgroup of finite index. 
Let $P \subset A$ be the inverse image of 
${\cal M}_{U,x}/{\cal O}^*_{U,x}$. 
Then $P$ is a local chart of ${\cal M}_{U}$ 
at a neighborhood of $x$ (cf.~\cite[(2.10)]{klog1}). 
It is clear that the morphism 
${\cal O}_{U,x}\otimes_{\mab Z}(P^{\rm gp}/Q^{\rm gp}) 
\owns a \otimes m \lom ad\log m \in  
\Lam^1_{U/V,x}$ is an isomorphism of 
${\cal O}_{U,x}$-modules. 
\par 
Set $W:=(\os{\circ}{V}
\times_{{\rm Spec}({\mab C}[Q])_{\rm an}}
{\rm Spec}({\mab C}[P])_{\rm an}, P^a)$. 
Then the natural morphism 
$W \lo V$ is log smooth 
by the implication (2)$\Lo$(1). 
Let $\pi \col U \lo W$ be also the natural morphism. 
We claim that the following morphism 
\begin{equation*} 
{\cal O}_{W,\os{\circ}{\pi}(x)}
\otimes_{\mab Z}(P^{\rm gp}/Q^{\rm gp}) 
\owns a \otimes m \lom ad\log m \in 
\Lam^1_{W/V,\os{\circ}{\pi}(x)} \quad 
(a \in {\cal O}_{W,\os{\circ}{\pi}(x)}, m \in P)
\tag{2.1.1}\label{eqn:ozpql}
\end{equation*} 
of ${\cal O}_{W,\os{\circ}{\pi}(x)}$-modules 
is an isomorphism. 
(For a homomorphism $R \lo S$ of commutative monoids 
with unit elements and for a point $y$ of 
${\rm Spec}({\mab C}[S])_{\rm an}$, 
I do not know whether 
$$\Lam^1_{({\rm Spec}({\mab C}[S])_{\rm an},S^a)/
({\rm Spec}({\mab C}[R])_{\rm an},R^a),y}\simeq 
{\cal O}_{{\rm Spec}({\mab C}[S])_{\rm an},y}
\otimes_{\mab Z}(S^{\rm gp}/{\rm Im}(R^{\rm gp}\lo S^{\rm gp}))$$ 
because I do not know whether there exists 
a derivation 
$${\cal O}_{{\rm Spec}({\mab C}[S])_{\rm an},y} 
\lo {\cal O}_{{\rm Spec}({\mab C}[S])_{\rm an},y}
\otimes_{\mab Z}
(S^{\rm gp}/{\rm Im}(R^{\rm gp}\lo S^{\rm gp}))$$ 
which extends the algebraic derivation 
$${\mab C}[S] \owns m \lom m\otimes m \in 
{\mab C}[S]\otimes_{\mab Z}
(S^{\rm gp}/{\rm Im}(R^{\rm gp}\lo S^{\rm gp})) \quad (m\in S).)$$
Set $P':={\mab N}^r\oplus Q$ and $W':=(\os{\circ}{V}
\times_{{\rm Spec}({\mab C}[Q])_{\rm an}}
{\rm Spec}({\mab C}[P'])_{\rm an}, P'{}^a)$. 
Then the natural morphism 
$P'{}^{\rm gp}\otimes_{\mab Z}{\mab Q} 
\lo P^{\rm gp}\otimes_{\mab Z}{\mab Q}$
is an isomorphism. 
Hence the natural morphism 
$p\col W \lo W'$ is log etale 
by the implication (2)$\Lo$(1). 
Furthermore, ${\rm Spec}({\mab C}[P'])_{\rm an}
=({\mab A}^1_{\mab C}, 
({\mab N}\owns 1 \lom z \in 
{\cal O}_{{\mab A}^1_{\mab C}})^a)^r\times
({\rm Spec}({\mab C}[Q])_{\rm an}, Q^a)$,  
where ${\mab A}^1_{\mab C}(={\mab C})$ 
is the line and $z$ is the holomorphic 
function ${\rm id} \col {\mab C} \lo {\mab C}$.  
Hence $\Lam^1_{W'/V,\os{\circ}{p}(\os{\circ}{\pi}(x))}
\simeq 
{\cal O}_{W',\os{\circ}{p}(\os{\circ}{\pi}(x))}
\otimes_{\mab Z}(P'{}^{\rm gp}/Q^{\rm gp})
={\cal O}_{W',\os{\circ}{p}(\os{\circ}{\pi}(x))}
\otimes_{\mab Z}(P^{\rm gp}/Q^{\rm gp})$. 
As in the algebraic case (\cite[(3.12)]{klog1}), 
for morphisms $g \col X \lo Y$ and $h \col Y \lo Z$ 
of fine log analytic spaces over ${\mab C}$, 
we have the following exact sequence 
\begin{equation*} 
g^*(\Lam^1_{Y/Z}) \lo \Lam^1_{X/Z} \lo 
\Lam^1_{X/Y} \lo 0.
\tag{2.1.2}\label{eqn:fstfes}
\end{equation*}
As in the algebraic case (\cite[(3.12)]{klog1}), 
as to the following conditions (i) and (ii), 
\par 
(i) $g$ is log smooth (resp.~log etale), 
\par 
(ii) $g^*(\Lam^1_{Y/Z})$ is 
a local direct factor of $\Lam^1_{X/Z}$ 
(resp.~$g^*(\Lam^1_{Y/Z}) \os{\sim}{\lo} \Lam^1_{X/Z}$), 
\parno
the following implications hold: 
(i)$\Lo$(ii); if $h\circ g$ is log smooth, 
then (ii)$\Lo$(i). Since $p$ is log etale, 
$\Lam^1_{W/V,\os{\circ}{\pi}(x)}=
p^*(\Lam^1_{W'/V})_{\os{\circ}{\pi}(x)}=
{\cal O}_{W,\os{\circ}{\pi}(x)}
\otimes_{\mab Z}(P^{\rm gp}/Q^{\rm gp})$. 
Therefore we have proved that 
the morphism (\ref{eqn:ozpql}) is an isomorphism. 
\par 
Now it is clear that 
$\Lam^1_{U/V,x}=\pi^*(\Lam^1_{W/V,\os{\circ}{\pi}(x)})$ 
and consequently the morphism $U \lo W$ is log etale by 
the implication (ii)$\Lo$(i). 
Since ${\cal M}_U$ is the pull-back of ${\cal M}_W$, 
the morphism 
$(\os{\circ}{U}, {\cal O}_U^*) 
\lo  
(\os{\circ}{W}, {\cal O}_W^*)$ is etale in our sense. 
The rest that we have to prove is 
that the morphism 
$\os{\circ}{U} \lo \os{\circ}{W}$
is locally isomorphic, which follows from 
the following lemma 
(\ref{lemm:tse}) (2) 
whose proof is not trivial.  
(Though the following lemma may be well-known 
(see [EGA IV-4] and \cite{sga1} for the algebraic case), 
we give the proof of it 
because we cannot find an appropriate 
reference in the analytic case.) 
\end{proof}

\begin{lemm}\label{lemm:tse} 
Let the notations be as above. 
Assume that the log analytic spaces 
$T_0$, $T$, $U$ and $V$ are trivial. 
Then the following hold$:$ 
\par 
$(1)$ Assume that 
$f$ is smooth$(:=$log smooth in our sense$)$. 
Then $f$ is flat. 
\par
$(2)$ The morphism 
$f$ is etale$(:=$log etale in our sense$)$ 
if and only if $f$ is locally isomorphic as 
a morphism of analytic spaces.  
\end{lemm} 
\begin{proof} 
(1):  (The proof of (1) is similar to that in 
[EGA IV-4, (17.5.1)].) 
By the proof of (\ref{prop:locdse}),  
$\Om^1_{U/V}$ is a locally  
projective ${\cal O}_U$-module of finite type. 
Let $x$ be a point of $U$ and 
let $y\in V$ be the image of $x$ by $f$. 
Let ${\mab C}\{u_1, \ldots, u_n\}$ be 
the ring of convergent series in 
$n$-variables over ${\mab C}$. 
Then it is well-known that 
${\mab C}\{u_1, \ldots, u_n\}$ $(n\in {\mab N})$ 
is a noetherian regular local ring. 
By the definition of an analytic space over ${\mab C}$, 
there exist nonnegative integers $m$ and $n$ such that 
${\cal O}_{U,x}={\mab C}\{u_1, \ldots, u_m\}/I_x$ 
and 
${\cal O}_{V,y}={\mab C}\{v_1, \ldots, v_n\}/I_y$ 
for some ideals $I_x$ and $I_y$ of 
${\mab C}\{u_1, \ldots, u_m\}$ and 
${\mab C}\{v_1, \ldots, v_n\}$, respectively. 
The morphism ${\cal O}_{V,y} \lo {\cal O}_{U,x}$ 
induces a surjection 
\begin{equation*} 
p\col {\mab C}\{u_1, \ldots, u_m, v_1,\ldots,v_n\}/I_y
{\mab C}\{u_1, \ldots, u_m, v_1,\ldots,v_n\} 
\lo  {\cal O}_{U,x}. 
\end{equation*} 
Let $I$ be the kernel of this surjection. 
Set 
$${\cal O}_{V,y}\{\ul{u}\} :=
{\mab C}\{u_1, \ldots, u_m, v_1,\ldots,v_n\}/I_y
{\mab C}\{u_1, \ldots, u_m, v_1,\ldots,v_n\}.$$  
Then  we claim that 
the natural morphism 
\begin{equation*} 
I/I^2\owns t \lom dt\otimes 1 \in  
\Om^1_{{\cal O}_{V,y}\{\ul{u}\}/{\cal O}_{V,y}}
\otimes_{{\cal O}_{V,y}\{\ul{u}\}}{\cal O}_{U,x} 
\tag{2.2.1}\label{eqn:iomo}
\end{equation*} 
has a left inverse. 
Indeed, the following exact sequence 
\begin{equation*} 
\begin{CD}
0 @>>> I/I^2 @>>> 
{\cal O}_{V,y}\{\ul{u}\}/I^2 
@>{p}>> {\cal O}_{U,x} @>>> 0 
\end{CD}
\end{equation*} 
of ${\cal O}_{V,y}$-modules 
is split since $U \lo V$ is smooth 
and since ${\cal O}_{V,y}\{\ul{u}\}/I^2$ defines 
a closed analytic space  
$V\times ({\mab A}^n_{\mab C})_{\rm an}$ 
in a neighborhood of $(x, O)$.  
Let $s \col {\cal O}_{U,x} \lo {\cal O}_{V,y}\{\ul{u}\}/I^2$
be a section of $p$.  
It is easy to check that 
$({\rm id}-s\circ p) \col {\cal O}_{V,y}\{\ul{u}\}/I^2 \lo 
I/I^2$ is a derivation over ${\cal O}_{V,y}$. 
Since $({\rm id}-s\circ p)\vert_{I/I^2}={\rm id}_{I/I^2}$,  
the derivation gives a desired left inverse of the morphism 
(\ref{eqn:iomo}). 
Therefore, by [EGA IV-1, Chapitre 0 (19.1.12)], 
there exist a generator $\{f_j\}_{j\in J}$ 
of finite elements of 
$I$ which generate $I/I^2$ 
and a subset $\{u_j\}_{j\in J} \subset \{u_1, \ldots, u_m\}$  
such that 
${\rm det}(\part f_i/\part u_j)_{i,j\in J} \notin 
{\mathfrak m}$, 
where ${\mathfrak m}$ is the maximal ideal of 
${\cal O}_{V,y}\{\ul{u}\}$.   
\par 
By noting that the ring 
${\mab C}\{u_1, \ldots, u_m\}$ 
is a noetherian regular local ring, 
the rest of the proof of (1) is 
similar to that of the implication a)$\Lo$b) in 
[EGA IV-4, (17.5.1)]: we directly use 
\cite[Chapitre 0 (10.2.4)]{ega31} instead of 
[EGA IV-3, (11.3.8)] used in [EGA IV-4, (17.5.1)]. 
\par 
(2): $\Lo)$: By the standard deformation theory 
(cf.~\cite[III (5.1)]{sga1}), 
we have $\Om^1_{U/V}=0$. 
As in [EGA IV-4 (17.4.1)] or \cite[I (3.5)]{mi}, the condition  
$\Om^1_{U/V}=0$ implies that $f$ is unramified 
in the classical sense. 
Consequently $f$ is etale in the classical sense 
by (1).  
Hence $f$ is locally isomorphic 
by \cite{car} mentioned in \cite[XII (3.3) a)]{sga1}. 
\par
The converse implication of (2) is obvious 
since $T_0=T$ as topological spaces. 
\end{proof}  

\begin{rema}
As in the algebraic case \cite[(3.6)]{klog1}, 
we can replace (a) in (\ref{prop:locdse}) 
by the following condition $(a)'$:
\par 
$(a)'$ the morphism $Q^{\rm gp}\lo P^{\rm gp}$ is injective 
$($resp.~and the morphism $Q^{\rm gp}\otimes_{\mab Z}{\mab Q} 
\lo P^{\rm gp}\otimes_{\mab Z}{\mab Q}$ is surjective$)$. 
\end{rema}

\par
As in \cite[(2.1.2) (iii)]{nale}, 
we can define a {\it Kummer morphism} 
of fs log analytic spaces over ${\mab C}$. 
The composite morphisms of two Kummer morphisms is 
of Kummer type. 

\begin{prop}[{cf.~\cite[(1.2)]{vimle}}]\label{prop:ked}
Let $f \col U \lo V$ be a morphism of 
fine log analytic spaces over ${\mab C}$. 
Then the following conditions $(1)$ and $(2)$ 
are equivalent$:$
\par 
$(1)$ $f$ is log etale and of Kummer type. 
\par 
$(2)$ There exists a local chart 
$Q \lo P$ of $f$ 
satisfying the following conditions 
$(a)_1$, $(a)_2$ and $(b)$ in 
$(\ref{prop:locdse}):$
\par 
$(a)_1$ $P$ and $Q$ are saturated.   
\par 
$(a)_2$ The morphism $Q \lo P$ is
injective and there exists 
a positive integer $n$ such that 
$P^n \subset {\rm Im}(Q \lo P)$. 
\end{prop}
\begin{proof}
By the proof of (\ref{prop:locdse}) and 
by the analytic analogue of 
\cite[(2.10)]{klog1}, 
we easily obtain (\ref{prop:ked}). 
\end{proof} 

As in \cite[(4.6)]{klog1}, we can define an exact morphism 
of fs log analytic spaces over ${\mab C}$. 
The composite morphism of two exact morphisms is exact. 
One can check by elementary calculations that 
the exactness is stable under the base change 
in the category of fine log 
analytic spaces over ${\mab C}$ and 
in that of fs log analytic spaces over ${\mab C}$.

\begin{coro}\label{cor:cbk}
For a morphism $f \col U \lo V$ of fs log analytic spaces 
over ${\mab C}$, $f$ is log etale and of Kummer type 
if and only if $f$ is log etale and exact. 
Consequently, the log etaleness of Kummer type  
is stable under the base change 
in the category of fs log 
analytic spaces over ${\mab C}$. 
\end{coro}
\begin{proof} 
(\ref{cor:cbk}) immediately follows from 
(\ref{prop:ked}). 
\end{proof}

\par
Let $Y$ be an fs log analytic space over ${\mab C}$. 
By (\ref{cor:cbk}), 
we obtain a topos $\wt{Y}^{\log}_{\rm et}$ 
which is the obvious analogue of 
the log etale topos of an fs log scheme (\cite[(2.2)]{nale}). 
In particular, for the trivial log analytic space 
$(\os{\circ}{Y}, {\cal O}_Y^*)$, 
we obtain the topos $\wt{\os{\circ}{Y}}_{\rm et}$. 
We call $\wt{Y}^{\log}_{\rm et}$ the 
({\it analytic}) {\it log etale topos} of $Y$. 
As in \cite[(2.5)]{nale}, we obtain 
the notion of the log geometric point 
of the topos $\wt{Y}^{\log}_{\rm et}$. 
\par
By using the local description of 
a Kummer log etale morphism $f \col U \lo V$ of 
fs log analytic spaces over ${\mab C}$ ((\ref{prop:ked}))
and using \cite[(1.3) (3)]{kn} and \cite[(1.2.1.1)]{kn}, 
the associated morphism 
$f^{\log} \col U^{\log} \lo V^{\log}$ 
is a local isomorphism of topological spaces 
by the same proof as that of \cite[(2.2)]{kn}. 
Hence we have a natural morphism  
\begin{equation*}
\be  \col \wt{Y}^{\log} \lo 
\wt{Y}^{\log}_{\rm et} 
\tag{2.5.1}\label{eqn:bedf}
\end{equation*}
of topoi. We often denote $\be$ by $\be_Y$. 
We also have a natural morphism  
\begin{equation*} 
\eps_{\rm an} \col 
\wt{Y}^{\log}_{\rm et} \lo \wt{\os{\circ}{Y}}_{\rm et} 
\tag{2.5.2}\label{eqn:yly}
\end{equation*}
of topoi. We sometimes denote 
$\eps_{\rm an}$ by $\eps_{Y}$. 
We have the following commutative diagram 
\begin{equation*}
\begin{CD}
\wt{Y}^{\log} @>{\be_Y}>> \wt{Y}^{\log}_{\rm et} \\
@V{\eps_{\rm cl}}VV @VV{\eps_{\rm an}}V \\
\wt{\os{\circ}{Y}} @>{\be_{\os{\circ}{Y}}}>> 
\wt{\os{\circ}{Y}}_{\rm et}. 
\end{CD}
\tag{2.5.3}\label{eqn:compmcz} 
\end{equation*}   
The direct image $\be_{\os{\circ}{Y}*}$ 
is an exact functor 
from the category of 
abelian sheaves in $\wt{\os{\circ}{Y}}$ 
to the category of abelian sheaves in 
$\wt{\os{\circ}{Y}}_{\rm et}$ by \cite[III (3.3)]{mi}
and commutes with the tensor product over ${\mab Z}$. 
It is trivial to check that  
$R\be_{\os{\circ}{Y}*}\be_{\os{\circ}{Y}}^{-1} 
\os{\sim}{\longleftarrow}{\rm id}$. 
More generally, 
\begin{equation*}
\be_{\os{\circ}{Y}*}\be_{\os{\circ}{Y}}^{-1} ={\rm id} 
\quad {\rm and} \quad 
\be_{\os{\circ}{Y}}^{-1}\be_{\os{\circ}{Y}*} ={\rm id} 
\tag{2.5.4}\label{cd:ssid} 
\end{equation*}
on the category 
of sheaves of sets in $\wt{\os{\circ}{Y}}_{\rm et}$ 
and $\wt{\os{\circ}{Y}}$, respectively.
Let ${\cal M}_{Y,\log}$ be a sheaf of monoids in 
$\wt{Y}^{\log}_{\rm et}$ which is associated to 
the presheaf $U \lom \Gam(U, {\cal M}_U)$ 
$(U\in Y^{\log}_{\rm et})$. 
Let ${\cal O}_{Y,\log}$ be 
the structure sheaf in 
$\wt{Y}^{\log}_{\rm et}$. 
Then we have a natural commutative diagram 
\begin{equation*} 
\begin{CD}
\eps^{-1}_{\rm an}\be_{\os{\circ}{Y}*}({\cal M}_{Y}) 
@>>> {\cal M}_{Y,\log} \\ 
@VVV @VVV \\ 
\eps^{-1}_{\rm an}\be_{\os{\circ}{Y}*}({\cal O}_{Y}) 
@>>> {\cal O}_{Y,\log}. 
\end{CD}
\tag{2.5.5}\label{cd:mclmlo}  
\end{equation*} 
Denote $({\cal M}_{Y,\log})^{\rm gp}$ simply by 
${\cal M}^{\rm gp}_{Y,\log}$.

\section{Proof of (\ref{theo:logint}) (1)}\label{sec:knctsc}
In this section we give the proof 
of (\ref{theo:logint}) (1). 
\par 
Let the notations be as in \S\ref{sec:lgelk}. 
The following is an analogue of \cite[(2.3)]{kn}:

\begin{lemm}[Analytic log Kummer sequence]
\label{lemm:logkum} 
For a positive integer $m$, the following sequence  
\begin{equation*} 
0 \lo  ({\mab Z}/m)(1) \lo {\cal M}^{\rm gp}_{Y,\log} 
\os{m\times}{\lo} 
{\cal M}^{\rm gp}_{Y,\log}  \lo  0   
\tag{3.1.1}\label{eqn:kumclet}
\end{equation*}  
is  exact in $\wt{Y}^{\log}_{\rm et}$. 
\end{lemm} 
\begin{proof} 
The obvious analytic analogue of 
the proof of \cite[(2.3)]{kn} works. 
\end{proof}
\par 
The natural morphism 
$\eps^{-1}_{\rm an}\be_{\os{\circ}{Y}*}({\cal M}_{Y}) \lo 
{\cal M}_{Y, \log}$ induces a morphism 
\begin{equation*}
\eps^{-1}_{\rm cl}({\cal M}_{Y}) \lo 
\be^{-1}_Y({\cal M}_{Y, \log}) 
\tag{3.1.2}\label{eqn:zmbe} 
\end{equation*} 
by (\ref{eqn:compmcz}) and (\ref{cd:ssid}).

\begin{prop}\label{prop:lclkum} 
Let $m$ be a positive integer. 
Let $E$ be an $m$-torsion abelian sheaf 
in $\wt{\os{\circ}{Y}}$. 
Then the canonical morphism 
\begin{equation*}
\bigwedge^k
({\cal M}^{\rm gp}_{Y}/{\cal O}_{Y}^*)
\otimes_{\mab Z}({\mab Z}/m)(-k)\otimes_{\mab Z}E {\lo} 
R^k\eps_{{\rm cl}*}(\eps^{-1}_{\rm cl}(E)) 
\quad (k \in {\mab Z}_{\geq 0}) 
\tag{3.2.1}\label{eqn:clkumz}
\end{equation*} 
by the use of {\rm (\ref{eqn:kumclet})} and {\rm (\ref{eqn:zmbe})} 
as in the paragraph between 
{\rm \cite[(2.3)]{kn} and \cite[(2.4)]{kn}}
is an isomorphism.  
\end{prop}
\begin{proof}
The proof is the same as that of \cite[(1.5)]{kn}. 
\end{proof}

\begin{prop}\label{prop:alclkumn} 
Let $m$ be a positive integer. 
Let $E$ be an $m$-torsion abelian sheaf 
in $\wt{\os{\circ}{Y}}_{\rm et}$. 
Then there exists a canonical isomorphism 
\begin{equation*}
\be_{\os{\circ}{Y}*}\{\bigwedge^k
({\cal M}^{\rm gp}_{Y}/{\cal O}_{Y}^*)\}
\otimes_{\mab Z}({\mab Z}/m)(-k)\otimes_{\mab Z}E
\os{\sim}{\lo} 
R^k\eps_{{\rm an}*}(\eps^{-1}_{\rm an}(E)) 
\quad (k \in {\mab Z}_{\geq 0}). 
\tag{3.3.1}\label{eqn:aclkumzn}
\end{equation*}
\end{prop} 
\begin{proof} 
As we said before, we have  
the notion of the log geometric point 
of the topos $\wt{Y}^{\log}_{\rm et}$. 
By using the analytic log Kummer sequence 
(\ref{eqn:kumclet}), 
the obvious analytic analogue of the proof of 
\cite[(2.4)]{kn} works. 
\end{proof}

\par 
We say that a log analytic space $U$ over ${\mab C}$ is 
{\it affine} 
if $\os{\circ}{U}$ is isomorphic to a closed analytic space 
of a polydisk over ${\mab C}$. 
\par 
The following is an analytic analogue of 
an algebraic constructible sheaf 
in a log etale topos in \cite[(3.3) 8]{nale}.
 
\begin{defi}\label{defi:ancon}
(1) Let $f \col Z \lo Y$ be a morphism of fs log analytic spaces 
over ${\mab C}$. 
We say that $f$ is {\it of log finite type} 
if, for any affine open log analytic space $U$ of $Y$ 
which has a global chart $P$ of ${\cal M}_U$,  
$f^{-1}(U)$ is the union of a finite number of 
open log analytic spaces $V_i$'s of $Z$ 
which have global charts $P\lo Q_i$ of 
$f\vert_{V_i} \col V_i \lo U$.    
\par 
(2) Let $A$ be a commutative ring with unit element. 
Let $K$ be a sheaf of $A$-modules in $\wt{Y}^{\log}_{\rm et}$. 
We say that $K$ is {\it constructible} 
if, for any affine open log analytic subspace $U$ 
of $Y$, there exist objects $V$ and $W$ of 
$U^{\log}_{\rm et}$ such that the structural morphisms 
$V \lo U$ and $W \lo U$ are of log finite type and 
such that $K\vert_U$ is isomorphic to 
the cokernel of a morphism 
$A_{V,U} \lo A_{W,U}$, where $A_{V,U}:=j_{V!}(A_V)$ 
and $A_{W,U}:=j_{W!}(A_W)$ for the structural morphisms 
$j_V \col V \lo U$ and $j_W \col W \lo U$, respectively. 
\end{defi}

\begin{defi}\label{defi:alelb}
Let $K$ be an abelian sheaf in $\wt{Y}^{\log}_{\rm et}$. 
We say that $K$ is {\it classical} 
(resp.~{\it classically constructible}) if  
there exists an abelian sheaf 
(resp.~constructible abelian sheaf) 
$E$ in $\wt{\os{\circ}{Y}}_{\rm et}$ such that 
$K$ is isomorphic to $\eps_{Y}^{-1}(E)$. 
We say that $K$ is {\it locally classical} 
(resp.~{\it locally classically constructible}) if, 
for any point $y$ of $\os{\circ}{Y}$, there exists 
an object $U \lo Y$ of $Y^{\log}_{\rm et}$ 
whose image in $Y$ contains $y$ and if  
$K\vert_U$ is 
classical (resp.~classically constructible). 
Let ${\rm D}^+_{{\rm lcl}{\textrm -}{\rm tor}}
(Y^{\log}_{\rm et})$  
(resp.~${\rm D}^+_{{\rm lclc}{\textrm -}{\rm tor}}
(Y^{\log}_{\rm et})$)
be the derived category of bounded below complexes of 
abelian sheaves in $\wt{Y}^{\log}_{\rm et}$ 
whose cohomology sheaves are 
locally classical and torsion 
(resp.~locally classically constructible and torsion). 
We obtain the similar notions in the algebraic case. 
\end{defi}

\begin{prop}\label{prop:cslc}
Let $A$ be a commutative ring with unit element. 
Let $K$ be a constructible sheaf of $A$-modules 
in $\wt{Y}^{\log}_{\rm et}$. Then $K$ is 
locally classically constructible. 
In particular, 
${\rm D}^+_{{\rm c}{\textrm -}{\rm tor}}(Y^{\log}_{\rm et})
\subset 
{\rm D}^+_{{\rm lclc}{\textrm -}{\rm tor}}
(Y^{\log}_{\rm et})$. 
The obvious algebraic analogue also holds. 
\end{prop}
\begin{proof} 
By using (\ref{prop:ked}), 
the proof is the same as that 
of \cite[(2.5.2)]{kn}.
\end{proof}

\begin{theo}\label{theo:kwle}  
The adjunction morphism 
\begin{equation*}
K^{\bul} \lo R\be_{Y*}(\be^{-1}_Y(K^{\bul})) 
\tag{3.7.1}\label{eqn:kbrbe}
\end{equation*}
for an object $K^{\bul}$ in 
${\rm D}^+_{{\rm lcl}{\textrm -}{\rm tor}}
(Y^{\log}_{\rm et})$
is an isomorphism. Consequently 
{\rm (\ref{theo:logint}) (1)} holds. 
\end{theo}
\begin{proof} 
We may assume that $K^{\bul}$ is 
a locally classical torsion abelian sheaf $K$ 
in $\wt{Y}^{\log}_{\rm et}$.  
By the lemma (\ref{lemm:lim}) (3) below, 
we may assume that $K$ is killed by 
a positive integer $m$.
Because the question is local, 
we may furthermore 
assume that $K=\eps_{\rm an}^{-1}(E)$ 
for some $m$-torsion abelian sheaf $E$ in 
$\wt{\os{\circ}{Y}}_{\rm et}$. 
By (\ref{prop:alclkumn}) and (\ref{prop:lclkum}), 
we see that  
the adjunction morphism 
${\rm id} \lo R\be_{Y*}\be^{-1}_Y$ induces an isomorphism 
\begin{equation*} 
R\eps_{{\rm an}*}(\eps_{\rm an}^{-1}(E)) 
\os{\sim}{\lo} 
R\be_{\os{\circ}{Y}*}R\eps_{{\rm cl}*}
(\eps_{\rm cl}^{-1}\be_{\os{\circ}{Y}}^{-1}(E)).
\tag{3.7.2}\label{eqn:nice}
\end{equation*} 
The same argument as that in \cite[p.~172, l.~1--7]{kn} 
tells us that   
(\ref{eqn:nice}) shows (\ref{theo:kwle}). 
\end{proof}

The following 
(3) is not included 
in a general theorem \cite[VI (5.1)]{sga4-2} 
since $\wt{Y}^{\log}_{\rm et}$ is not algebraic in 
general.

\begin{lemm}\label{lemm:lim} 
$(1)$ Let $\{F_{\lam}\}_{\lam\in \Lam}$ be an inductive 
system of abelian sheaves in $\wt{Y}^{\log}$. Then 
\begin{equation*} 
\vil_{\lam \in \Lam}R\eps_{{\rm cl}*}(F_{\lam})
=
R\eps_{{\rm cl}*}(\vil_{\lam \in \Lam}F_{\lam}). 
\tag{3.8.1}\label{eqn:eplm}
\end{equation*}
\par
$(2)$ Let $\{F_{\lam}\}_{\lam\in \Lam}$ be an inductive 
system of abelian sheaves in $\wt{Y}^{\log}_{\rm et}$. 
Then 
\begin{equation*} 
\vil_{\lam \in \Lam}R\eps_{{\rm an}*}(F_{\lam})
=
R\eps_{{\rm an}*}(\vil_{\lam \in \Lam}F_{\lam}). 
\tag{3.8.2}\label{eqn:anlime}
\end{equation*}
\par 
$(3)$ Let $\{F_{\lam}\}_{\lam\in \Lam}$ be an inductive 
system of abelian sheaves in $\wt{Y}^{\log}$. Then 
\begin{equation*} 
\vil_{\lam \in \Lam}R\be_{Y*}(F_{\lam})
=
R\be_{Y*}(\vil_{\lam \in \Lam}F_{\lam}). 
\tag{3.8.3}\label{eqn:anbelm}
\end{equation*}
\end{lemm}
\begin{proof} 
(1): By the proper base change theorem for 
locally compact spaces (\cite[II (4.11.1)]{go}), 
we may assume that $\os{\circ}{Y}$ is a point. 
Then $Y^{\log}=({\mab S}^1)^r$ for some $r \in {\mab N}$ 
(\cite[(1.3)]{kn}). 
Since ${\mab S}^1$ is compact, 
the cohomology of an abelian sheaf 
with compact support on $Y^{\log}$ 
and the usual cohomology of an abelian sheaf on 
$Y^{\log}$ are the same. 
Because to take the direct limit of abelian sheaves and 
to take the cohomology 
with compact support are commutative 
(\cite[II (4.11.2)]{go}, \cite[III (5.1)]{iv}), we obtain (1). 
\par
(2): The question is local on $\os{\circ}{Y}_{\rm et}$. 
Hence we may assume that $Y$ has a global chart 
$P \lo {\cal O}_Y$. In this case, we can calculate 
the higher direct image $R^h\eps_{{\rm an}*}$ 
$(h \in {\mab N})$ by the sheafied version 
of the cohomology of 
the group ${\rm Hom}(P^{\rm gp}, \hat{\mab Z}(1))$ 
as in the algebraic case in \cite[(4.7.1)]{nale}. 
Hence (2) follows.
\par 
(3): As in the argument in 
\cite[p.~172, l.~1--7]{kn}, we have only to prove that    
\begin{equation*} 
R\eps_{{\rm an}*}\vil_{\lam \in \Lam}R\be_{Y*}(F_{\lam})
=
R\eps_{{\rm an}*}R\be_{Y*}(\vil_{\lam \in \Lam}F_{\lam}). 
\tag{3.8.4}\label{eqn:rvilm}
\end{equation*}
By (1) and (2) and by the quasi-equivalence of 
$\be_{\os{\circ}{Y}*}$, 
the left hand side on (\ref{eqn:rvilm}) 
is equal to 
\begin{align*} 
\vil_{\lam \in \Lam}
R\eps_{{\rm an}*}R\be_{Y*}(F_{\lam})  & =
\vil_{\lam \in \Lam}
R\be_{\os{\circ}{Y}*}R\eps_{{\rm cl}*}(F_{\lam}) \\
{} & = R\be_{\os{\circ}{Y}*}R\eps_{{\rm cl}*}
(\vil_{\lam \in \Lam}F_{\lam}) \\
{} & =
R\eps_{{\rm an}*}R\be_{Y*}(\vil_{\lam \in \Lam}F_{\lam}). 
\end{align*}
\end{proof}

\par
Because I do not know whether a log exponential 
sequence in $\wt{Y}^{\log}_{\rm et}$ exists 
when ${\cal M}_Y$ is nontrivial, 
I do not know the answer of the following problem: 

\begin{prob}\label{prob:isont}
Let ${\rm D}^+_{\rm lcl}(Y^{\log}_{\rm et})$ 
be the derived category of bounded below complexes of 
abelian sheaves in $\wt{Y}^{\log}_{\rm et}$ 
whose cohomology sheaves are 
locally classical. Is the adjunction morphism 
(\ref{eqn:kbrbe}) isomorphic 
for an object $K^{\bul}$ in 
${\rm D}^+_{\rm lcl}(Y^{\log}_{\rm et})$?
\end{prob}
 
\section{Proof of (\ref{theo:logint}) (2)}\label{sec:bcg}
In this section 
we give a base change theorem (\ref{theo:basec}) 
below which is a generalization of 
(\ref{theo:logint}) (2). 
\par 
Let $X$ be an fs log scheme over ${\mab C}$ 
whose underlying scheme is locally of finite type 
over ${\mab C}$. 
Let ${\cal O}_X$ 
and ${\cal O}_{X_{\rm an}}$ 
be the structure sheaves in $\wt{\os{\circ}{X}}_{\rm et}$ 
and $\wt{\os{\circ}{X}}_{\rm an}$, respectively. 
Let ${\cal M}_X$ be the log structure in 
$\wt{\os{\circ}{X}}_{\rm et}$. 
Let $\eta \col \wt{\os{\circ}{X}}_{\rm an} 
\lo \wt{\os{\circ}{X}}_{\rm et}$    
be the natural morphism of topoi. 
Let $\eta^{*}({\cal M}_{X})$ be the log structure 
on $\os{\circ}{X}_{\rm an}$ which is associated to 
the composite morphism 
$\eta^{-1}({\cal M}_{X}) \lo 
\eta^{-1}({\cal O}_X) \lo 
{\cal O}_{X_{\rm an}}$.  
As in \cite{kn}, 
we call the fs log analytic space 
$X_{\rm an}:=
(\os{\circ}{X}_{\rm an}, \eta^*({\cal M}_X))$ 
the {\it log analytic space associated to} $X$.  
We call  $\eta^{*}({\cal M}_{X})(=:{\cal M}_{X_{\rm an}})$  
the {\it analytification} 
of ${\cal M}_{X}$. 
Also, as in \cite{kn}, $X^{\log}_{\rm an}$ 
denotes the real blow up of $X_{\rm an}$. 
\par
Let $\eta^{\log} \col 
\wt{X}^{\log}_{\rm an} 
\lo \wt{X}^{\log}_{\rm et}$   
be the natural morphism of topoi (\cite[(2.1), (2.2)]{kn}). 
Then we have the following commutative diagram of topoi 
(\cite[p.~171]{kn}):  
\begin{equation*}
\begin{CD}
\wt{X}^{\log}_{\rm an}  @>{\eta^{\log}}>> 
\wt{X}^{\log}_{\rm et} \\
@V{\eps_{\rm cl}}VV @VV{\eps_{\rm et}}V \\
\wt{\os{\circ}{X}}_{\rm an}  @>{\eta}>> 
\wt{\os{\circ}{X}}_{\rm et}. 
\end{CD} 
\tag{4.0.1}\label{eqn:loganet}
\end{equation*} 
We sometimes denote $\eps_{\rm et}$ by $\eps_{X}$, 
and $\eta^{\log}$ and $\eta$ by 
$\eta^{\log}_X$ and $\eta_X$, respectively.

\begin{lemm}\label{lemm:contf} 
Let $f \col U\lo V$ be a morphism of fs log schemes 
over ${\mab C}$ whose underlying schemes are locally of 
finite type over ${\mab C}$. 
Let $f_{\rm an} \col U_{\rm an} \lo V_{\rm an}$ 
be the associated morphism of fs log analytic spaces 
over ${\mab C}$.  
If $f$ is log etale and of Kummer type, then 
$f_{\rm an}$ is so.  
\end{lemm} 
\begin{proof}
By the local descriptions of 
algebraic and analytic log etale morphisms of Kummer type 
(\cite[(3.5)]{klog1} 
(cf.~\cite[p.~369]{nale}, \cite[(1.2)]{vimle}), 
(\ref{prop:ked})),   
(\ref{lemm:contf}) 
immediately follows.
\end{proof}

\par 
Because a functor $U \lom U_{\rm an}$ 
$(U \in {X}^{\log}_{\rm et})$ defines  
a continuous functor 
$X^{\log}_{\rm et} \lo (X_{\rm an})^{\log}_{\rm et}$ 
by (\ref{lemm:contf}), 
we have a morphism 
\begin{equation*}
\eta_{\rm et} \col \wt{(X_{\rm an})}{}^{\log}_{\rm et} \lo 
\wt{X}^{\log}_{\rm et} 
\tag{4.1.1}\label{eqn:uxet}
\end{equation*}
of topoi. 
We often denote $\eta_{\rm et}$ by $\eta_{X,{\rm et}}$.  
Then we have the following commutative diagram: 

\begin{equation*} 
\begin{CD}
\wt{X}^{\log}_{\rm an} @>{\be_{X_{\rm an}}}>> 
\wt{(X_{\rm an})}{}^{\log}_{\rm et} @>{\eta_{X,{\rm et}}}>> 
\wt{X}^{\log}_{\rm et} \\ 
@V{\eps_{\rm cl}}VV @V{\eps_{\rm an}}VV 
@V{\eps_{\rm et}}VV \\  
\wt{\os{\circ}{X}}_{\rm an} 
@>{\be_{\os{\circ}{X}_{\rm an}}}>> 
\wt{(\os{\circ}{X}_{\rm an})}_{\rm et} 
@>{\eta_{\os{\circ}{X},{\rm et}}}>> 
\wt{\os{\circ}{X}}_{\rm et}.
\end{CD}
\tag{4.1.2}\label{cd:bgc}
\end{equation*}
By the definition of $\eta^{\log}$ and $\eta$, 
we have two equalities:
\begin{equation*}
\eta^{\log}=  \eta_{X,{\rm et}} \circ \be_{X_{\rm an}}, 
\quad 
\eta=  \eta_{\os{\circ}{X},{\rm et}} 
\circ \be_{\os{\circ}{X}_{\rm an}}.
\tag{4.1.3}\label{eqn:eeb}
\end{equation*}

\begin{prop}\label{prop:qqcs}
Let $A$ be a commutative ring with unit element. 
Assume that $\os{\circ}{X}$ is quasi-compact 
and quasi-separated. 
Then,  
for a constructible sheaf $K$ of $A$-modules 
in $\wt{X}^{\log}_{\rm et}$, 
$\eta_{X,{\rm et}}^{-1}(K)$ is a constructible 
sheaf of $A$-modules in 
$\wt{(X_{\rm an})}{}^{\log}_{\rm et}$. 
Consequently $\eta_{X,{\rm et}}^{-1}$ induces a functor 
${\rm D}^+_{{\rm c}{\textrm -}{\rm tor}}(X^{\log}_{\rm et}) 
\lo 
{\rm D}^+_{{\rm c}{\textrm -}{\rm tor}}
((X_{\rm an})^{\log}_{\rm et})$. 
\end{prop}
\begin{proof}
The proof is easy. 
\end{proof}

For a log scheme $X$ and a log analytic space $Y$ over 
${\mab C}$, 
${\rm D}^+(X^{\log}_{\rm et})$ and 
${\rm D}^+(Y^{\log}_{\rm et})$ denote
the derived categories of bounded below complexes of 
abelian sheaves in $\wt{X}^{\log}_{\rm et}$ and 
$\wt{Y}^{\log}_{\rm et}$, respectively.

\begin{lemm}\label{lemm:gzw}
Let ${g} \col {Z} \lo {W}$ be a morphism of 
analytic spaces over ${\mab C}$. 
As in {\rm \S2}, $\wt{{Z}}$ and $\wt{{W}}$ denote 
$\wt{{Z}}_{\rm cl}$ and $\wt{{W}}_{\rm cl}$, 
respectively.  
Let $\be_Z \col \wt{{Z}} \lo \wt{{Z}}_{\rm et}$ 
and $\be_W \col \wt{{W}} \lo \wt{{W}}_{\rm et}$ 
be the morphisms of topoi defined in 
{\rm (\ref{eqn:bedf})} for the trivial log analytic spaces 
$Z$ and $W$, respectively. 
Let ${g}_{\rm cl} \col \wt{{Z}} \lo \wt{{W}}$
and 
${g} \col 
\wt{{Z}}_{\rm et} \lo 
\wt{{W}}_{\rm et}$ be the induced morphisms of topoi  
by ${g} \col {Z} \lo {W}$. 
Let $G^{\bul}$ be an object of ${\rm D}^+(Z_{\rm et})$. 
Then the following base change morphism 
\begin{equation*} 
\be_W^{-1}R{g}_{*}(G^{\bul}) \lo 
R{g}_{{\rm cl}*}(\be_Z^{-1}(G^{\bul})) 
\tag{4.3.1}\label{eqn:bgg}
\end{equation*} 
is an isomorphism.
\end{lemm}
\begin{proof}
We may assume that 
$G^{\bul}$ is an abelian sheaf $G$ in $\wt{Z}_{\rm et}$. 
Let $G' \lo H'$ be a morphism of abelian sheaves in 
$\wt{{W}}$. Then this morphism is an isomorphism 
if and only if so is the induced morphism 
$\be_{W*}(G') \lo \be_{W*}(H')$. 
Hence it suffices to prove that 
the morphism 
\begin{equation*} 
\be_{W*}\be_W^{-1}R{g}_{*}(G) 
\lo 
\be_{W*}Rg_{{\rm cl}*}(\be_Z^{-1}(G))
\tag{4.3.2}\label{eqn:gzst}
\end{equation*}  
is an isomorphism.  
The left hand side on (\ref{eqn:gzst}) 
is equal to $R{g}_{*}(G)$.
On the other hand, the target of 
the morphism (\ref{eqn:gzst}) is equal to 
\begin{equation*}
R\be_{W*}R{g}_{{\rm cl}*}(\be_Z^{-1}(G))
=R{g}_{*}R\be_{Z*}(\be_Z^{-1}(G))=R{g}_{*}(G).
\end{equation*}
Now we complete the proof. 
\end{proof}

The following is a variant of 
Artin-Grothendieck's' base change theorem 
(\cite[XVI (4.1)]{sga4-3}).

\begin{coro}\label{coro:vag}
Let $g \col Z \lo {W}$ 
be a morphism of schemes which 
are locally of finite type over ${\mab C}$. 
Let $\eta_{Z,{\rm et}} \col \wt{({Z}_{\rm an})}_{\rm et} 
\lo \wt{{Z}}_{\rm et}$ 
and 
$\eta_{W,{\rm et}} 
\col \wt{({W}_{\rm an})}_{\rm et} \lo \wt{{W}}_{\rm et}$ 
be the morphisms of topoi defined in {\rm (\ref{eqn:uxet})} 
for the trivial log schemes 
$Z$ and $W$, respectively. 
Let ${g}_{\rm an} \col 
\wt{{(Z_{\rm an})}}_{\rm et} \lo 
\wt{{(W_{\rm an})}}_{\rm et}$ 
and 
${g} \col \wt{Z}_{\rm et} \lo \wt{W}_{\rm et}$ 
be the induced morphisms  
of topoi by ${g} \col {Z} \lo {W}$. 
Assume that ${g}$ is of finite type. 
Let $K^{\bul}$ be 
an object of 
${\rm D}^+_{{\rm c}{\textrm -}{\rm tor}}(Z_{\rm et})$. 
Then the following base change morphism 
\begin{equation*} 
\eta_{{W},{\rm et}}^{-1}R{g}_*(K^{\bul}) 
\lo 
R{g}_{{\rm an}*}(\eta^{-1}_{{Z},{\rm et}}(K^{\bul})). 
\tag{4.4.1}\label{eqn:bctm}
\end{equation*} 
is an isomorphism. 
\end{coro}
\begin{proof}
We may assume that 
$K^{\bul}$ is a 
constructible torsion abelian sheaf $K$ in 
$\wt{Z}_{\rm et}$. 
Consider the following commutative diagram$:$ 
\begin{equation*} 
\begin{CD}
\wt{{Z}}_{\rm an} @>{\be_{{Z}_{\rm an}}}>> 
\wt{({Z}_{\rm an})}_{\rm et} 
@>\eta_{{Z},{\rm et}}>>
\wt{{Z}}_{\rm et}\\ 
@V{({g}_{\rm an})_{\rm cl}}VV 
@V{{g}_{\rm an}}VV @V{g}VV \\
 \wt{{W}}_{\rm an}
@>{\be_{{W}_{\rm an}}}>>
\wt{({W}_{\rm an})}_{\rm et} 
@>\eta_{{W},{\rm et}}>> \wt{{W}}_{\rm et}. 
\end{CD}
\tag{4.4.2}\label{cd:bcxty}
\end{equation*}  
Since $\be_{W_{\rm an}}^{-1}$ is exact, it suffices to 
prove that the morphism 
\begin{equation*} 
\be_{W_{\rm an}}^{-1}
\eta_{{W},{\rm et}}^{-1}R{g}_*(K) 
\lo 
\be_{W_{\rm an}}^{-1}
R{g}_{{\rm an}*}(\eta^{-1}_{{Z},{\rm et}}(K)) 
\tag{4.4.3}\label{eqn:pbbctm}
\end{equation*} 
is an isomorphism. 
By (\ref{lemm:gzw}) we see 
that the target of the morphism 
(\ref {eqn:pbbctm}) 
is equal to 
$R{(g_{\rm an})}_{{\rm cl}*}(
(\eta_{{Z},{\rm et}}\be_{Z_{\rm an}})^{-1}(K))$. 
By Artin-Grothendieck's' base change theorem 
(\cite[XVI (4.1)]{sga4-3}), this complex is equal to 
$(\eta_{{W},{\rm et}}\be_{W_{\rm an}})^{-1}R{g}_{*}(K)$. 
Now it is clear that the morphism 
(\ref {eqn:pbbctm}) 
is an isomorphism. 
\end{proof}

\begin{lemm}\label{lemm:xye}
Let $f \col X \lo Y$ be a morphism 
of fs log schemes over ${\mab C}$. 
Assume that $f$ is strict. 
Let $E^{\bul}$ be an object of 
${\rm D}^+(\wt{\os{\circ}{X}}_{\rm et})$. 
Then the following base change morphism 
\begin{equation*}
\eps_Y^{-1}R\os{\circ}{f}_*(E^{\bul}) 
\lo Rf_*(\eps_X^{-1}(E^{\bul}))
\tag{4.5.1}\label{eqn:esrf}
\end{equation*} 
is an isomorphism.
The obvious analogue 
for the analytic log etale topoi of  
fs log analytic spaces over ${\mab C}$ also holds. 
\end{lemm}
\begin{proof} 
We may assume that 
$E^{\bul}$ is an abelian sheaf $E$ in 
$\wt{\os{\circ}{X}}_{\rm et}$. 
Because the question is local on $Y$,  
we may assume that $Y$ has a global chart $P$, where $P$ 
is an fs monoid. 
Since $f$ is strict, $X$ also 
has a global chart $P$. Set 
$I(X):={\rm Hom}_{\rm gp}(P^{\rm gp}, \wh{\mab Z}(1))$ 
as in \cite[p.~376]{nale}. 
Then, by \cite[(4.6.1)]{nale}, 
the category of abelian sheaves in $\wt{X}^{\log}_{\rm et}$ 
is equivalent to the category of $I(X)$-modules in 
$\wt{\os{\circ}{X}}_{\rm et}$. 
The sheaf $\eps^{-1}_X(E)$ corresponds to $E$ 
with trivial $I(X)$-action (\cite[(4.7) (i)]{nale}). 
Let $\eps^{-1}_X(E) \lo I^{\bul}$ be 
an injective resolution of $\eps^{-1}_X(E)$.  
Because it is clear that $I^{\bul}$ can be considered as 
an injective resolution of $E$, 
the morphism (\ref{eqn:esrf}) is an isomorphism 
by \cite[(4.7) (ii)]{nale}. 
\par
The proof for the analytic case is the same. 
\end{proof}


\par 
The following is a main result of this section. 

\begin{theo}[Base change theorem]\label{theo:basec} 
Let $f \col X \lo Y$ be a morphism 
of fs log schemes over ${\mab C}$ 
whose underlying schemes are locally of finite type over 
${\mab C}$. 
Assume that $\os{\circ}{f}$ is of finite type. 
Let $f_{\rm an} \col \wt{(X_{\rm an})}{}^{\log}_{\rm et} \lo 
\wt{(Y_{\rm an})}{}^{\log}_{\rm et}$ 
be the associated morphism to $f$. 
Consider the following commutative diagram$:$ 
\begin{equation*} 
\begin{CD}
\wt{(X_{\rm an})}{}^{\log}_{\rm et} 
@>{\eta_{X,{\rm et}}}>> \wt{X}^{\log}_{\rm et}
 \\ 
@V{f_{\rm an}}VV @VV{f}V \\
\wt{(Y_{\rm an})}{}^{\log}_{\rm et} @>{\eta_{Y,{\rm et}}}>> 
\wt{Y}^{\log}_{\rm et}. 
\end{CD}
\tag{4.6.1}\label{cd:bcxy}
\end{equation*}  
Let $K^{\bul}$ be an object of 
${\rm D}^+_{{\rm c}{\textrm -}{\rm tor}}(X^{\log}_{\rm et})$. 
If $f$ is log injective 
{\rm (\cite[(5.5.1)]{nale})}, 
then the following base change morphism 
\begin{equation*} 
\eta_{Y,{\rm et}}^{-1}Rf_*(K^{\bul}) \lo 
Rf_{{\rm an}*}(\eta^{-1}_{X,{\rm et}}(K^{\bul})) 
\tag{4.6.2}\label{eqn:bcm}
\end{equation*} 
is an isomorphism. 
\end{theo} 
\begin{proof} 
We may assume that $K^{\bul}$ is 
a constructible torsion abelian sheaf 
$K$ in $\wt{X}^{\log}_{\rm et}$. 
We decompose the morphism 
$f$ by the composite morphism $X \lo 
(\os{\circ}{X},\os{\circ}{f}{}^*({\cal M}_Y)) \lo Y$. 
Set $Z:=(\os{\circ}{X},\os{\circ}{f}{}^*({\cal M}_Y))$ and 
let $g \col Z \lo Y$ be the morphism above. 
\par 
First we claim that 
the base change morphism  (\ref{eqn:bcm}) for $g$ is 
an isomorphism. Indeed, 
let $\{Z_i\}_{i \in I}$ be a 
Kummer log etale covering of $Z$ 
such that  the pull-backs of $K$ to 
$Z_i$ for all $i\in I$ 
are classically constructible and torsion
((\ref{prop:cslc})). 
Set $Z_0:=\coprod_{i\in I}Z_i$ and 
$Z_{\bul}:={\rm cosk}_0^Z(Z_0)$ $(\bul \in {\mab N})$. 
By the cohomological descent 
in \cite[${\rm V}^{\rm bis}$]{sga4-2} 
and by the same proof as that of 
the crystalline base change theorem  
\cite[Chapitre V Th\'{e}or\`{e}me 3.5.1]{bb} 
(in our proof, we do not need to assume that 
$f$ is quasi-separated nor that the index set $I$ is finite 
because we do not need to consider 
the derived functor $\ul{\ul{\mab L}}$ in \cite{bb}), 
we have only to prove that the base change morphism 
(\ref{eqn:bcm}) for each morphism $Z_n \lo Y$ $(n\in {\mab N})$ 
is an isomorphism. 
(Here we have used the Convention (3) implicitly 
in the cohomological descent.)   
Consequently we may assume that 
$K=\eps^{-1}_{Z}(G)$ for a 
constructible torsion abelian sheaf 
$G$ in $\wt{\os{\circ}{Z}}_{\rm et}$. 
Then we have the following formulae 
\begin{align*} 
\eta_{Y,{\rm et}}^{-1}Rg_*(\eps_{Z}^{-1}(G)) & 
\os{(\ref{lemm:xye})}{=} 
\eta_{Y,{\rm et}}^{-1}\eps_{Y}^{-1}R\os{\circ}{g}_*(G)= 
\eps_{Y_{\rm an}}^{-1}\eta_{\os{\circ}{Y},{\rm et}}^{-1}
R\os{\circ}{g}_*(G) \tag{4.6.3}\label{eqn:agbf}\\ 
{} & 
\os{(\ref{coro:vag})}{=}
\eps_{Y_{\rm an}}^{-1}R\os{\circ}{g}_{{\rm an}*}
(\eta_{\os{\circ}{Z},{\rm et}}^{-1}(G)) 
\os{(\ref{lemm:xye})}{=} 
Rg_{{\rm an}*}(\eps_{Z_{\rm an}}^{-1}
\eta_{\os{\circ}{Z},{\rm et}}^{-1}(G)) \\
{} & =
Rg_{{\rm an}*}(\eta_{Z,{\rm et}}^{-1}\eps_{Z}^{-1}(G)). 
\end{align*}
Here the numbers above equalities mean that the equalities 
follow from the statements numbered.  
Hence we have proved the claim.
\par
Let $h \col X \lo Z$ be the morphism above.
Secondly, we claim that 
the base change morphism  (\ref{eqn:bcm}) for $h$ is 
an isomorphism. In fact we claim that 
the morphism (\ref{eqn:bcm}) for $h$ is an isomorphism for 
a bounded below complex of 
abelian sheaves in $\wt{X}^{\log}_{\rm et}$. 
Indeed, let $z$ be a closed point of 
$\os{\circ}{Z}(=\os{\circ}{X})$. 
Because the problem is local on $Z$, 
we may assume that there exists a global chart 
$P_Z \lo P_X$ of $h$, where $P_Z$ and $P_X$ are 
fs monoids.
Set $I_X:={\rm Hom}(P_X^{\rm gp}, \hat{\mab Z}(1))$, 
$I_Z:={\rm Hom}(P_Z^{\rm gp}, \hat{\mab Z}(1))$ and 
$I_{h}:={\rm Ker}(I_X \lo I_Z)$. 
Let $K$ be an abelian sheaf in $\wt{X}^{\log}_{\rm et}$; 
$K$ corresponds to an abelian sheaf $\wt{K}$ in 
$\wt{\os{\circ}{X}}_{\rm et}$ with continuous 
$I_X$-action (\cite[(4.6.1)]{nale}). 
Then, by 
\cite[(4.7.1)]{nale}, 
it is easy to check that 
$R^qh_*(K)_{z(\log)}=
{\rm Map}_{c, I_X/I_h}(I_Z, H^q(I_h, \wt{K}_z))$.  
On the other hand, by the analytic analogues of 
\cite[(4.6.1)]{nale} and \cite[(4.7.1)]{nale}, 
$R^qh_{{\rm an}*}(K)_{z_{\rm an}(\log)}
={\rm Map}_{c, I_X/I_h}(I_Z,
H^q(I_h,\eta^{-1}_{X,{\rm et}}(\wt{K})_{z_{\rm an}}))
={\rm Map}_{c, I_X/I_h}(I_Z,H^q(I_h,\wt{K}_{z}))$, 
where $z_{\rm an}$ is the point of $\os{\circ}{Z}_{\rm an}$ 
corresponding to $z$. 
\par 
Finally (\ref{theo:basec}) follows from 
\cite[(5.5.2)]{nale}. 
Indeed, because $f$ is log injective, 
so is $h$. 
Hence [loc.~cit.] tells us that 
$R^qh_*(K)$ is constructible, which permits us to use 
the base change theorem (\ref{theo:basec}) 
for $g$ as follows:
\begin{align*} 
Rf_{{\rm an}*}(\eta^{-1}_{X,{\rm et}}(K))
& =Rg_{{\rm an}*}Rh_{{\rm an}*}
(\eta^{-1}_{X,{\rm et}}(K))
=Rg_{{\rm an}*}\eta^{-1}_{Z,{\rm et}}Rh_*(K) 
\tag{4.6.4}\label{eqn:cbcm}\\ 
{} & =\eta^{-1}_{Y,{\rm et}}Rg_*Rh_*(K)
=\eta^{-1}_{Y,{\rm et}}Rf_*(K).
\end{align*}
\end{proof}

\begin{coro}\label{theo:up}
Let $K^{\bul}$ be an object of 
${\rm D}^+_{{\rm c}{\textrm -}{\rm tor}}(X^{\log}_{\rm et})$. 
Then the adjunction morphism 
$$K^{\bul} \lo 
R\eta_{{\rm et}*}(\eta^{-1}_{\rm et}(K^{\bul}))$$ 
is an isomorphism. Consequently 
{\rm (\ref{theo:logint}) (2)} holds.  
\end{coro}
\begin{proof}
We may assume that $K^{\bul}$ is 
a constructible torsion abelian sheaf $K$ 
in $\wt{X}^{\log}_{\rm et}$. 
We may assume that $\os{\circ}{X}$ 
is quasi-compact. 
By (\ref{theo:basec}), 
for any object $U \in X^{\log}_{\rm et}$, 
$H^p(\wt{(U_{\rm an})}{}^{\log}_{\rm et}, 
\eta^{-1}_{\rm et}(K\vert_U))=
H^p(\wt{U}^{\log}_{\rm et}, K\vert_U)$ $(p\in {\mab N})$. 
Hence $\eta_{{\rm et}*}(\eta^{-1}_{\rm et}(K))=K$ and 
$R^p\eta_{{\rm et}*}(\eta^{-1}_{\rm et}(K))=0$ for 
$p>0$ (cf.~\cite[III (2.11) (a)]{mi}). 
\end{proof}

\begin{coro}[{\cite[(2.6)]{kn}}]\label{coro:reo}
Let $K^{\bul}$ be an object of 
${\rm D}^+_{{\rm c}{\textrm -}{\rm tor}}(X^{\log}_{\rm et})$.  
Then the adjunction morphism 
$$K^{\bul} \lo R\eta^{\log}_{*}(\eta^{\log,-1}(K^{\bul}))$$ 
is an isomorphism. Consequently 
{\rm (\ref{theo:cole})} holds. 
\end{coro}
\begin{proof} 
(\ref{coro:reo}) immediately follows 
from the following equality 
\begin{equation*}
\eta^{\log}= \eta_{\rm et}\circ \be,   
\tag{4.8.1}\label{eqn:elbce}
\end{equation*} 
and from (\ref{prop:qqcs}), (\ref{prop:cslc}), 
(\ref{theo:kwle}) and 
(\ref{theo:up}).
\end{proof}

Though we do not use the log proper base change theorem 
(\cite[(5.1)]{nale}) 
nor the following analytic log proper base change theorem 
in the proof of (\ref{theo:basec}), 
we may use them 
(however we shall not use them in this paper).

\begin{theo}[Analytic log proper base change theorem]
\label{alpbt}
Let 
\begin{equation*} 
\begin{CD} 
X_4 @>{g'}>> X_1 \\ 
@V{f'}VV @VV{f}V \\ 
X_3 @>{g}>> X_2
\end{CD}
\end{equation*} 
be a cartesian diagram of fs log analytic spaces over 
${\mab C}$. 
Assume that $\os{\circ}{f}$  is proper and 
that the condition in 
{\rm \cite[(5.1)]{nale}} holds, that is, 
for any two points $x_1\in  \os{\circ}{X}_1$ and 
$x_3\in  \os{\circ}{X}_3$ lying over a point 
$x_2 \in  \os{\circ}{X}_2$, and for $P_i:=
({\cal M}_{X_i}/{\cal O}^*_{X_i})_{x_i}$, 
the inverse image of $P_1\oplus P_3$ 
by the morphism 
$P^{\rm gp}_2\owns a \lom (a, a^{-1}) \in  
P^{\rm gp}_1\oplus P^{\rm gp}_3$
is $\{1\}$. Then, for a bounded below complex $K^{\bul}$ 
of abelian sheaves in $(\wt{X}_1)^{\log}_{\rm et}$ 
whose cohomology sheaves are torsion, 
the base change morphism 
$g^{-1}Rf_*(K^{\bul}) \lo Rf'_*(g'{}^{-1}(K^{\bul}))$ 
is an isomorphism. 
\end{theo}
\begin{proof} 
We point out only the different points from 
the proof of \cite[(5.1)]{nale}. 
We may assume that $K^{\bul}$ is a torsion abelian sheaf $K$ 
in $(\wt{X}_1)^{\log}_{\rm et}$.  
\par 
We consider the case corresponding to 
\cite[Claim (5.2.1)]{nale}: assume that 
$g$ is strict. 
Let $x_2$ be a point of $\os{\circ}{X}_2$. 
We may assume that $X_2$ has a global chart 
$P\lo {\cal O}_{X_2}$, where $P$ is an fs monoid. 
By (\ref{prop:ked}), 
$R^qf_*(K)_{x_2(\log)}=
\vil_{n,U}H^q((X_1\times_{X_2}U)
\times_{{\rm Spec}({\mab C}[P])_{\rm an}}
{\rm Spec}({\mab C}[P^{1/n}])_{\rm an}, K)$, 
where $\os{\circ}{U}$ is an open neighborhood of $x_2$. 
As in \cite[(5.2.1)]{nale}, we may assume that 
the log structure of $X_2$ is trivial. 
Then, by the classical proper base change theorem 
(\cite[II (4.11.1)]{go}), we may assume that 
$\os{\circ}{f}$ is the identity. 
Then, as in the proof of 
(\ref{theo:basec}), two direct calculations 
by using the analytic analogues of \cite[(4.6.1)]{nale} 
and [loc.~cit.~(4.7.1)], we obtain (\ref{alpbt}) 
when $g$ is strict. 
\par
The rest of the proof is the same as that of 
\cite[(5.1)]{nale} 
by replacing ${\rm Spec}(?)$ in [loc.~cit.] 
by ${\rm Spec}(?)_{\rm an}$. 
\end{proof}

\begin{rema}\label{rema:lslt} 
(1) Let $g \col Z \lo W$ be a morphism of 
fs log analytic spaces over ${\mab C}$. 
Consider the following commutative diagram 
\begin{equation*} 
\begin{CD}
\wt{Z}^{\log} @>{\beta_Z}>> \wt{Z}^{\log}_{\rm et} \\ 
@V{g^{\log}}VV @VV{g}V \\
\wt{W}^{\log} @>{\beta_W}>> 
\wt{W}^{\log}_{\rm et}. 
\end{CD}
\end{equation*}  
Then the example \cite[(2.7)]{kn} tells us that 
the base change morphism 
\begin{equation*} 
\beta_{W}^{-1}Rg_*(K) \lo 
Rg^{\log}_{*}(\beta^{-1}_{Z}(K)) 
\tag{4.10.1}\label{eqn:esbcm}
\end{equation*}
for a constructible torsion abelian sheaf 
$K$ in $\wt{Z}^{\log}_{\rm et}$ 
is not an isomorphism in general. 
\par 
Let $f \col X \lo Y$ be a morphism of fs log schemes  
whose underlying schemes are locally of 
finite type. The difference between 
the higher direct images of 
$f^{\log}_{\rm an} \col 
\wt{X}^{\log}_{\rm an} \lo 
\wt{Y}^{\log}_{\rm an}$ and 
$f \col \wt{X}^{\log}_{\rm et}\lo \wt{Y}^{\log}_{\rm et}$  
pointed out in [loc.~cit.] 
arises from the difference between those of 
$f^{\log}_{\rm an} \col 
\wt{X}^{\log}_{\rm an} \lo 
\wt{Y}^{\log}_{\rm an}$ and 
$f_{\rm an} \col 
\wt{(X_{\rm an})}{}^{\log}_{\rm et} \lo 
\wt{(Y_{\rm an})}{}^{\log}_{\rm et}$.   
\par
(2) C.~Nakayama has informed me that 
he and T.~Kajiwara have proved that 
the base change morphism in 
\cite[(2.7.1)]{kn} is an isomorphism 
if $f$ is log injective. 
It is reasonable to expect 
that the obvious analytic analogue of 
this also holds. 
\par 
(3) I expect that the base change morphism 
(\ref{eqn:bcm}) is an isomorphism even if 
$f$ is not necessarily log injective.  
\par 
(4) Because we introduce the topos 
$\wt{Y}^{\log}_{\rm et}$ in \S\ref{sec:lgelk}, 
we can give a proof of 
\cite[(2.6)]{kn} in (\ref{coro:reo}) by using   
Artin-Grothendieck's' base change theorem 
(\cite[XVI (4.1)]{sga4-3}) in a full form 
and only by using the analytic log Kummer sequence; 
we have not used 
the algebraic log Kummer sequence nor 
the log exponential sequence in \cite{kn} 
to obtain \cite[(2.6)]{kn}, 
though we have used the same proof 
as that of \cite[(1.5)]{kn} 
in the proof of (\ref{prop:lclkum}). 
\end{rema}

\par 
We conclude this section by 
giving a logarithmic version of 
Grauert-Remmert's theorem 
(\cite{grat}, \cite[XI (4.3)]{sga4-3}), 
which is of independent interest. 
\par
Let $f \col X \lo Y$ be a morphism of fine log schemes 
(resp.~fine log analytic spaces)
over ${\mab C}$. 
We say that $f$ is {\it finite} 
if the morphism $\os{\circ}{f} 
 \col \os{\circ}{X} \lo \os{\circ}{Y}$ is finite.  
\par
Let $M$ be a finite abelian group. 
Let $X$ be an fs log scheme over ${\mab C}$ 
whose underlying scheme is locally of finite type 
over ${\mab C}$. Then, by (\ref{theo:logint}) (2), 
we have 
\begin{equation*}
H^1(X^{\log}_{\rm et}, M) =
H^1((X_{\rm an})^{\log}_{\rm et}, M), 
\tag{4.10.2}\label{eqn:mxx}
\end{equation*}
which makes us expect the following 
((\ref{eqn:mxx}) immediately follows 
from the following):

\begin{theo}[Log Grauert-Remmert's theorem]\label{theo:lgr}
Let $X$ be a fine log scheme over ${\mab C}$ 
whose underlying scheme is 
locally of finite type over ${\mab C}$. 
Let ${\rm Fet}(X)$ be the category of 
finite log etale coverings over $X$ 
and ${\rm Fet}(X_{\rm an})$ the similar category for 
$X_{\rm an}$. Then the functor of analytifications 
\begin{equation*} 
{\rm Fet}(X)\owns X' \lom  X'_{\rm an} \in  
{\rm Fet}(X_{\rm an}) 
\tag{4.11.1}\label{eqn:nxaxy}
\end{equation*}
gives an equivalence of categories.
\par 
Assume furthermore that ${\cal M}_X$ is saturated. 
Let ${\rm Fket}(X)$ be the category of 
finite log etale coverings of Kummer types over $X$ 
and ${\rm Fket}(X_{\rm an})$ the similar category for 
$X_{\rm an}$. 
Then the functor of analytifications 
\begin{equation*} 
{\rm Fket}(X)\owns X' \lom  X'_{\rm an} \in  
{\rm Fket}(X_{\rm an}) 
\tag{4.11.2}\label{eqn:knxeaxy}
\end{equation*}
gives an equivalence of categories. 
\end{theo}
\begin{proof}
Because the proof for 
(\ref{eqn:knxeaxy}) is similar to that for 
(\ref{eqn:nxaxy}), 
we give only the proof for (\ref{eqn:nxaxy}). 
\par  
Let $X'$ and $X''$ be two objects of 
${\rm Fet}(X)$. For simplicity of notation, 
set $Y:=X_{\rm an}$, $Y':=X'_{\rm an}$ and 
$Y'':=X''_{\rm an}$. Then it is easy to check 
that the natural map 
\begin{equation*} 
{\rm Hom}_X(X',X'') \lo {\rm Hom}_{Y}(Y',Y'')   
\tag{4.11.3}\label{eqn:nhmxy}
\end{equation*}
is injective. 
Let $h$ be an element of 
${\rm Hom}_Y(Y',Y'')$. 
By \cite[XI (4.3)]{sga4-3} the morphism  
$\os{\circ}{h} \col 
\os{\circ}{Y}{}' \lo \os{\circ}{Y}{}''$ 
is associated to a morphism 
$\os{\circ}{g} \col \os{\circ}{X}{}' \lo \os{\circ}{X}{}''$ 
of schemes over $\os{\circ}{X}$. 
Let $x\in \os{\circ}{X}{}'({\mab C})$ be a 
closed point of  $\os{\circ}{X}{}'$. 
Let $P' \lo {\cal M}_{X'}\vert_{\os{\circ}{U}{}'}$ and 
$P'' \lo {\cal M}_{X''}\vert_{\os{\circ}{U}{}''}$ be 
two charts, 
where $\os{\circ}{U}{}'$ and $\os{\circ}{U}{}''$ 
are etale neighborhoods 
of $\ol{x}$ and $\ol{\os{\circ}{g}(x)}$.
Then we have two composite 
surjective homomorphisms 
$$P'{}^{\rm gp} \lo  
{\cal M}^{\rm gp}_{X',\ol{x}}/{\cal O}^*_{X',\ol{x}}
\os{\sim}{\lo} 
{\cal M}^{\rm gp}_{Y',\ol{x}}/{\cal O}^*_{Y',\ol{x}}$$ 
and 
$$P''{}^{\rm gp} \lo  
{\cal M}^{\rm gp}_{X'',
\ol{\os{\circ}{g}(x)}}/{\cal O}^*_{X'',\ol{\os{\circ}{g}(x)}}
\os{\sim}{\lo} 
{\cal M}^{\rm gp}_{Y'',\ol{\os{\circ}{h}(x)}}/
{\cal O}^*_{Y'',\ol{\os{\circ}{h}(x)}}.$$ 
Let $K' \subset P'{}^{\rm gp}$ and 
$K'' \subset P''{}^{\rm gp}$ be 
the kernels of the homomorphisms above and 
let $P'_1 \subset P'{}^{\rm gp}/K'$ and 
$P''_1 \subset P''{}^{\rm gp}/K''$ 
be the inverse images of 
${\cal M}_{Y',\ol{x}}/{\cal O}^*_{Y',\ol{x}}$ 
and 
${\cal M}_{Y'',\ol{\os{\circ}{h}(x)}}/
{\cal O}^*_{Y'',\ol{\os{\circ}{h}(x)}}$, respectively.  
Then $P'_1$ and $P''_1$ are also the inverse images of 
${\cal M}_{X',\ol{x}}/{\cal O}^*_{X',\ol{x}}$ 
and 
${\cal M}_{X'',\ol{\os{\circ}{g}(x)}}
/{\cal O}^*_{X'',\ol{\os{\circ}{g}(x)}}$.
The morphism $h$ induces a homomorphism 
$P''_1 \lo P'_1$ and this gives a local chart of $h$ 
at $\ol{x}$ and $\ol{\os{\circ}{h}(x)}$.  
By \cite[(2.10)]{klog1}, 
$(P'_1)^a\vert_{\os{\circ}{U}{}'_1}=
{\cal M}_{Y'}\vert_{\os{\circ}{U}{}'_1}$ 
and $(P''_1)^a\vert_{\os{\circ}{U}{}''_1}
={\cal M}_{Y''}\vert_{\os{\circ}{U}{}''_1}$ 
for some etale neighborhoods $\os{\circ}{U}{}'_1$ and 
$\os{\circ}{U}{}''_1$ 
of $\ol{x}$ and $\ol{\os{\circ}{g}(x)}$, respectively,  
with a morphism 
$\os{\circ}{g}{}_1 \col \os{\circ}{U}{}'_1 \lo \os{\circ}{U}{}''_1$ 
over $\os{\circ}{g}$.
Hence we have a morphism 
$\os{\circ}{g}{}^{-1}_1({\cal M}_{Y''}\vert_{\os{\circ}{U}{}''_1}) 
\lo 
{\cal M}_{Y'}\vert_{\os{\circ}{U}{}'_1}$. 
Therefore the map (\ref{eqn:nhmxy}) is 
etale-locally surjective. 
In fact, it is surjective since it is injective. 
Now we have proved that the functor  
(\ref{eqn:nxaxy}) is fully faithful. 
\par 
The rest is to prove that the functor  
(\ref{eqn:nxaxy}) is essentially surjective. 
Let $Z$ be an object of ${\rm Fet}(Y)$.
Since the problem is local, we may 
assume that there exists a chart $P \lo Q$ 
of the morphism $Z \lo Y$ such that the induced morphism 
$P^{\rm gp}\otimes_{\mab Z}{\mab Q} 
\lo Q^{\rm gp}\otimes_{\mab Z}{\mab Q}$ is an isomorphism  
and we may assume that
the induced morphism 
$\os{\circ}{Z} \lo \os{\circ}{Y}
\times_{{\rm Spec}({\mab C}[P])_{\rm an}}
{\rm Spec}({\mab C}[Q])_{\rm an}$ is etale 
((\ref{prop:locdse})). 
Then this morphism is finite. Indeed, 
by the assumption, the composite morphism 
$\os{\circ}{Z} \lo \os{\circ}{Y}
\times_{{\rm Spec}({\mab C}[P])_{\rm an}}
{\rm Spec}({\mab C}[Q])_{\rm an} \lo \os{\circ}{Y}$ is finite.  
Furthermore, the morphism  
$\os{\circ}{Y}\times_{{\rm Spec}({\mab C}[P])_{\rm an}}
{\rm Spec}({\mab C}[Q])_{\rm an} \lo \os{\circ}{Y}$ 
is separated by the stability of the separation under 
the base change in the algebraic case and by 
the GAGA of the separation (\cite[XII (3.1)]{sga1}). 
Hence, by the obvious analytic analogue of 
\cite[(6.1.5)]{ega2}, 
we see that the morphism $\os{\circ}{Z} \lo \os{\circ}{Y}
\times_{{\rm Spec}({\mab C}[P])_{\rm an}}
{\rm Spec}({\mab C}[Q])_{\rm an}$ is finite. 
By (\ref{lemm:tse}) (2), $\os{\circ}{Z}$ is 
a finite covering space of 
$\os{\circ}{Y}\times_{{\rm Spec}({\mab C}[P])_{\rm an}}
{\rm Spec}({\mab C}[Q])_{\rm an}$ as topological spaces.
Now, by \cite[XI (4.3)]{sga4-3}, there exists a 
finite etale covering 
$\os{\circ}{X}{}' \lo \os{\circ}{X}
\times_{{\rm Spec}({\mab C}[P])}{\rm Spec}({\mab C}[Q])$ 
such that $\os{\circ}{X}{}'_{\rm an}=\os{\circ}{Z}$. 
Endow $\os{\circ}{X}{}'$ with the log structure 
associated to the morphism $Q \lo {\cal O}_{\os{\circ}{X}{}'}$ 
and let $X'$ be the resulting log scheme. 
Then $X'_{\rm an}=Z$.    
\end{proof}

\begin{rema}
Let $X$ be an fs log scheme over ${\mab C}$ 
whose underlying scheme is 
locally of finite type over ${\mab C}$.  
Though a finite abelian Galois covering 
of $X^{\log}_{\rm an}$ as topological spaces 
is obtained by a finite abelian Galois covering of $X$ by 
(\ref{theo:cole}), 
I know nothing about the nonabelian 
case. 
\end{rema}

\section{Analytic vs. algebraic log Kummer sequences}
\label{sec:knct}
The aim in this section 
is to give a commutative diagram 
(\ref{cd:kugcret}) below which 
compares the calculations of certain higher direct images 
by the use of 
the analytic and algebraic log Kummer sequences.  
\par 
Let $X$ be an fs log scheme 
over ${\mab C}$ whose underlying scheme $\os{\circ}{X}$ is 
locally of finite type over ${\mab C}$. 
\par 
Let ${\cal M}_{X,\log}$ be a sheaf of monoids 
which is associated to the presheaf 
$U \lom \Gam(U, {\cal M}_U)$ $(U \in X^{\log}_{\rm et})$.  
\par 
Let ${\cal O}_{X,\log}$ be the structure sheaf 
in $\wt{X}^{\log}_{\rm et}$. 
Then we have natural commutative diagrams  
\begin{equation*} 
\begin{CD}
\eps^{-1}_{\rm et}({\cal M}_X) 
@>>> {\cal M}_{X,\log} \\ 
@VVV @VVV \\ 
\eps^{-1}_{\rm et}({\cal O}_X) 
@>>> {\cal O}_{X,\log} 
\end{CD}
\tag{5.0.1}\label{cd:esmclmll}  
\end{equation*} 
and 
\begin{equation*} 
\begin{CD}
\eta^{-1}_{\rm et}({\cal M}_{X,\log}) 
@>>> {\cal M}_{X_{\rm an},\log} \\ 
@VVV @VVV \\ 
\eta^{-1}_{\rm et}({\cal O}_{X,\log}) 
@>>> {\cal O}_{X_{\rm an},\log}.  
\end{CD}
\tag{5.0.2}\label{cd:mclmloal}   
\end{equation*} 
Let $\eta^*_{\rm et}({\cal M}_{X,\log}) \in 
\wt{(X_{\rm an})}{}^{\log}_{\rm et}$ 
be the associated log structure to 
the composite morphism 
$\eta^{-1}_{\rm et}({\cal M}_{X,\log}) \lo 
\eta^{-1}_{\rm et}({\cal O}_{X,\log})  \lo 
{\cal O}_{X_{\rm an},\log}.$ 
Let 
\begin{equation*} 
\eta^{-1}_{\rm et}({\cal M}_{X,\log}) 
\lo 
\eta^*_{\rm et}({\cal M}_{X,\log})  
\tag{5.0.3}\label{eqn:e-1mle}
\end{equation*} 
be the natural morphism.  Set 
\begin{equation*}
\eta^{\log *}({\cal M}_{X,\log}):=
\be^{-1}_{X_{\rm an}}(\eta_{\rm et}^*({\cal M}_{X,\log})) 
\tag{5.0.4}\label{eqn:stea}
\end{equation*}
by abuse of notation.

\begin{defi}\label{defi:pbmx}
We call $\eta_{\rm et}^*({\cal M}_{X,\log})$ 
the {\it analytification} of ${\cal M}_{X,\log}$. 
\end{defi}
\parno
The upper horizontal morphism 
in (\ref{cd:esmclmll}) 
induces a morphism $
\eps_{\rm an}^{-1}\eta^{-1}_{\os{\circ}{X},{\rm et}}({\cal M}_X) 
\lo \eta^{-1}_{\rm et}({\cal M}_{X,\log})$. 
Composing this morphism with the morphism 
(\ref{eqn:e-1mle}), 
we have a morphism 
\begin{equation*}
\eps_{\rm an}^{-1}\eta^{-1}_{\os{\circ}{X},{\rm et}}
({\cal M}_X) 
\lo \eta^*_{\rm et}({\cal M}_{X,\log}). 
\tag{5.1.1}\label{eqn:epet*}
\end{equation*} 
In fact we have a natural morphism 
\begin{equation*}
\eps_{\rm an}^{-1}\eta^*_{\os{\circ}{X},{\rm et}}({\cal M}_X) 
\lo \eta^*_{\rm et}({\cal M}_{X,\log}). 
\tag{5.1.2}\label{eqn:*epet*}
\end{equation*} 
By the upper morphism of (\ref{cd:mclmloal}), we have 
a natural morphism 
\begin{equation*}
\eta^*_{\rm et}({\cal M}_{X,\log}) \lo 
{\cal M}_{X_{\rm an},\log}.  
\tag{5.1.3}\label{eqn:etmx}
\end{equation*} 
By the definition of $\eta^{\log*}({\cal M}_{X,\log})$,  
by (\ref{eqn:e-1mle}) and by (\ref{eqn:elbce}), 
we have a natural morphism 
\begin{equation*}
\eta^{\log, -1}({\cal M}_{X,\log}) 
\lo \eta^{\log*}({\cal M}_{X,\log}).   
\tag{5.1.4}\label{eqn:-1*et}
\end{equation*}
By the definition of $\eta^{\log*}({\cal M}_{X,\log})$,  
we have the following formula 
\begin{equation*} 
\eta^{\log *}({\cal M}_{X,\log})/\be^{-1}_{X_{\rm an}}
({\cal O}_{X_{\rm an},\log}^*)= 
\eta^{\log, -1}({\cal M}_{X,\log}/{\cal O}_{X,\log}^*).  
\tag{5.1.5}\label{eqn:etmoiv}
\end{equation*} 
Pulling back the morphism 
(\ref{eqn:epet*}) by the functor $\be^{-1}_{X_{\rm an}}$, 
we have a morphism 
\begin{equation*}
\eps_{\rm cl}^{-1}\eta^{-1}({\cal M}_X)  
\lo \eta^{\log *}({\cal M}_{X,\log}). 
\tag{5.1.6}\label{eqn:epeta}
\end{equation*} 
Let the notation be as in the beginning of 
\S\ref{sec:lgelk}. 
Let $\gam_{\os{\circ}{Y}} 
\col  
\wt{\os{\circ}{Y}}_{\rm et} 
\lo \wt{\os{\circ}{Y}}
(\not=\wt{\os{\circ}{Y}}_{\rm cl})$ 
be the natural morphism of topoi. 
Let ${\cal O}_{\os{\circ}{Y}}$ 
be the structure sheaf of 
the topos $\wt{\os{\circ}{Y}}
(\not=\wt{\os{\circ}{Y}}_{\rm cl})$. 
Because there exists a natural morphism 
$\eps^{-1}_{\rm an}\gam_{\os{\circ}{X}_{\rm an}}^{-1}
({\cal O}_{\os{\circ}{X}_{\rm an}})
\lo {\cal O}_{X_{\rm an},\log}$, 
the morphism (\ref{eqn:epeta}) induces a morphism 
\begin{equation*}
\eps_{\rm cl}^{-1}\eta^*({\cal M}_X)  
\lo \eta^{\log *}({\cal M}_{X,\log}). 
\tag{5.1.7}\label{eqn:stepet}
\end{equation*} 
Composing the morphism 
$\be^{-1}_{X_{\rm an}}((\ref{eqn:etmx}))$ with 
the morphism above, we obtain the following morphism 
\begin{equation*}
\eps_{\rm cl}^{-1}\eta^*({\cal M}_X)  
\lo \be^{-1}_{X_{\rm an}}({\cal M}_{X_{\rm an},\log}). 
\tag{5.1.8}\label{eqn:ster}
\end{equation*}

\par
Let $m$ be a positive integer. 
Recall the algebraic log Kummer sequence 
\begin{equation*}
0 \lo ({\mab Z}/m)(1) \lo 
{\cal M}_{X,\log}^{\rm gp} \os{m\times}{\lo} 
{\cal M}_{X,\log}^{\rm gp} \lo 0 
\tag{5.1.9}\label{eqn:logkums}
\end{equation*}
in $\wt{X}{}^{\log}_{\rm et}$ 
(\cite[(2.3)]{kn}). 
Then we have the following commutative diagram 
\begin{equation*}
\begin{CD} 
0 @>>> ({\mab Z}/m)(1) @>>> 
{\cal M}^{\rm gp}_{X_{\rm an},\log}  
@>{m\times}>> {\cal M}^{\rm gp}_{X_{\rm an},\log} @>>> 0 \\ 
@. @| @A{(\ref{eqn:etmx})^{\rm gp}}AA 
@AA{(\ref{eqn:etmx})^{\rm gp}}A  \\
0 @>>> ({\mab Z}/m)(1) @>>> 
\eta^{*}_{\rm et}({\cal M}^{\rm gp}_{X,\log})  
@>{m\times}>>  
\eta^{*}_{\rm et}({\cal M}^{\rm gp}_{X,\log}) @>>> 0 \\ 
@. @| @A{(\ref{eqn:e-1mle})^{\rm gp}}AA 
@AA{(\ref{eqn:e-1mle})^{\rm gp}}A \\
0 @>>> ({\mab Z}/m)(1) @>>> 
\eta^{-1}_{\rm et}({\cal M}^{\rm gp}_{X,\log})  
@>{m\times}>>  
\eta^{-1}_{\rm et}({\cal M}^{\rm gp}_{X,\log}) @>>> 0    
\end{CD}
\tag{5.1.10}\label{eqn:kumkum}
\end{equation*} 
of exact sequences (the exactness of the middle sequence 
is easy to check). 
\par 
Let $E$ be an $m$-torsion abelian sheaf 
in $\wt{\os{\circ}{X}}_{\rm et}$. 
Then Kato and Nakayama have proved that 
the log Kummer sequence 
(\ref{eqn:logkums}) gives the following canonical 
isomorphism (\cite[(2.4)]{kn}): 
\begin{equation*}
\bigwedge^k({\cal M}^{\rm gp}_{X}
/{\cal O}^*_{X})\otimes_{\mab Z}E(-k) 
\os{\sim}{\lo} R^k\eps_{{\rm et}*}(\eps^{-1}_{\rm et}(E)) 
\quad (k \in {\mab Z}_{\geq 0}).  
\tag{5.1.11}\label{cd:knis}
\end{equation*}
Here, in (\ref{cd:knis}),
we change the turn of 
the tensor product in [loc.~cit.]  
because the cup product is 
usually taken by the left cup product 
of a fundamental section 
(see \cite{rz} for example). 
\par 
Let $K^{\bul}$ be an object of 
${\rm D}^+(X^{\log}_{\rm et})$.
Then we have the following base change morphism 
\begin{equation*} 
\eta^{-1}R\eps_{{\rm et}*}(K^{\bul})
\lo R\eps_{{\rm cl}*}(\eta^{\log,-1}(K^{\bul})). 
\tag{5.1.12}\label{eqn:regrtf}
\end{equation*} 
In particular, for a nonnegative integer $k$, 
we have the following morphism 
\begin{equation*} 
\eta^{-1}R^k\eps_{{\rm et}*}(K^{\bul})
\lo R^k\eps_{{\rm cl}*}(\eta^{\log,-1}(K^{\bul})). 
\tag{5.1.13}\label{eqn:regkrtf}
\end{equation*} 
Hence we have a canonical morphism 
\begin{equation*} 
R^k\eps_{{\rm et}*}(K^{\bul})
\lo R\eta_*(R^k\eps_{{\rm cl}*}(\eta^{\log,-1}(K^{\bul}))).
\tag{5.1.14}\label{eqn:grtf}
\end{equation*} 

\begin{prop}\label{theo:keylelb} 
Let $m$ be a positive integer and 
let $E$ be an $m$-torsion abelian sheaf 
in $\wt{\os{\circ}{X}}_{\rm et}$. 
Let $k$ be a nonnegative integer. 
Then there exists 
the following commutative diagram 
\begin{equation*} 
\begin{CD} 
R\eta_*(\bigwedge^k(\eta^*({\cal M}^{\rm gp}_X) 
/{\cal O}^*_{X_{\rm an}})\otimes_{\mab Z}\eta^{-1}(E)(-k))
@>{\us{\sim}{R\eta_*((\ref{eqn:clkumz}))}}>> 
R\eta_*(R^k\eps_{{\rm cl}*}(\eps_{\rm cl}^{-1}\eta^{-1}(E))) \\
@A{}AA @AA{(\ref{eqn:grtf})}A \\
\bigwedge^k({\cal M}^{\rm gp}_{X}
/{\cal O}^*_{X})\otimes_{\mab Z}E(-k) 
@>{\us{\sim}{(\ref{cd:knis})}}>> 
R^k\eps_{{\rm et}*}(\eps^{-1}_{\rm et}(E)),  
\end{CD}
\tag{5.2.1}\label{cd:kugcret}
\end{equation*} 
where the left vertical morphism above is 
induced by the adjunction morphism 
${\rm id} \lo R\eta_*\eta^{-1}$. 
Furthermore, if $E$ is constructible, then 
the left vertical morphism is an isomorphism. 
\end{prop} 
\begin{proof} 
As to the commutativity of the diagram (\ref{cd:kugcret}), 
it suffices to prove that 
the following diagram is commutative: 
\begin{equation*} 
\begin{CD} 
\bigwedge^k(\eta^*({\cal M}^{\rm gp}_X) 
/{\cal O}^*_{X_{\rm an}})\otimes_{\mab Z}\eta^{-1}(E)(-k))
@>{\us{\sim}{(\ref{eqn:clkumz})}}>> 
R^k\eps_{{\rm cl}*}(\eps_{\rm cl}^{-1}\eta^{-1}(E)) \\
@| @AA{(\ref{eqn:grtf})}A \\
\eta^{-1}(\bigwedge^k({\cal M}^{\rm gp}_{X}
/{\cal O}^*_{X})\otimes_{\mab Z}E(-k)) 
@>{\us{\sim}{\eta^{-1}((\ref{cd:knis}))}}>> 
\eta^{-1}R^k\eps_{{\rm et}*}(\eps^{-1}_{\rm et}(E)),  
\end{CD}
\tag{5.2.2}\label{cd:kpbcret}
\end{equation*} 
Using (\ref{eqn:kumkum}), 
we have the following commutative diagram of triangles 
\begin{equation*}
\begin{CD} 
R\eps_{{\rm cl}*}(({\mab Z}/m)(1)) @>>> 
R\eps_{{\rm cl}*}(\be^{-1}({\cal M}^{\rm gp}_{X_{\rm an},\log}))  
@>{m\times}>> 
R\eps_{{\rm cl}*}(\be^{-1}({\cal M}^{\rm gp}_{X_{\rm an},\log})) 
@>{+1}>>  \\ 
@| @AAA @AAA  \\ 
R\eps_{{\rm cl}*}(({\mab Z}/m)(1)) @>>> 
R\eps_{{\rm cl}*}(\eta^{\log,-1}({\cal M}^{\rm gp}_{X,\log}))  
@>{m\times}>>  
R\eps_{{\rm cl}*}(\eta^{\log,-1}({\cal M}^{\rm gp}_{X,\log})) 
@>{+1}>>  \\ 
@AAA @AAA @AAA \\
\eta^{-1}R\eps_{{\rm et}*}(({\mab Z}/m)(1)) @>>> 
\eta^{-1}R\eps_{{\rm et}*}({\cal M}^{\rm gp}_{X,\log})  
@>{m\times}>>  
\eta^{-1}R\eps_{{\rm et}*}({\cal M}^{\rm gp}_{X,\log}) @>{+1}>>    
\end{CD}
\tag{5.2.3}\label{eqn:kumbk}
\end{equation*} 
(Here we have used the Convention (2).)
In particular, we have the commutativity of the diagram 
(\ref{cd:kpbcret}) for the case $k=1$ and $E={\mab Z}/m$. 
We leave the reader to the detail of 
the rest of the proof of the commutativity of the diagram 
(\ref{cd:kpbcret}) because 
it is a routine work by using the Godement resolution 
of an abelian sheaf in a topos with enough points and 
using the definition of the cup product. 
\par 
Assume now that $E$ is constructible. 
Since 
$({\cal M}^{\rm gp}_{X}/{\cal O}^*_{X})
\otimes_{\mab Z}{\mab Z}/m$ 
is a constructible torsion abelian sheaf in 
$\wt{\os{\circ}{X}}_{\rm et}$, 
the left vertical morphism in (\ref{cd:kugcret}) is 
an isomorphism by Artin-Grothendieck's comparison theorem 
(\cite[XVI (4.1)]{sga4-3}) as used in 
the proof of \cite[(2.6)]{kn}. 
\end{proof}

\begin{rema}
Using some results and some arguments 
in \cite[(2.6)]{kn},  using (\ref{theo:keylelb}) 
but without using 
the log exponential sequence in \cite{kn}, 
we can give a proof of [loc.~cit.] again. 
But we omit it because it resembles 
the proof in [loc.~cit.]. 
\end{rema}

Let 
\begin{equation*} 
R^k\eps_{{\rm et}*}(K^{\bul})
\lo R\eta_{\os{\circ}{X},{\rm et}*}
(R^k\eps_{{\rm an}*}(\eta^{-1}_{X,{\rm et}}(K^{\bul})))
\tag{5.3.1}\label{eqn:agnrtf}
\end{equation*} 
be an analogous morphism to (\ref{eqn:grtf}). 
\begin{prop}
Let $m$ be a positive integer and 
let $E$ be an $m$-torsion abelian sheaf 
in $\wt{\os{\circ}{X}}_{\rm et}$. 
Let $k$ be a nonnegative integer. 
Then the following diagram is commutative: 
\begin{equation*} 
\begin{CD} 
R\eta_{\os{\circ}{X},{\rm et}*}
(\bigwedge^k(\eta^{-1}_{\os{\circ}{X},{\rm et}*}
({\cal M}^{\rm gp}_X/{\cal O}^*_{X})
\otimes_{\mab Z}\eta^{-1}_{\os{\circ}{X},{\rm et}}(E)(-k))
@>{\us{\sim}{R\eta_{\os{\circ}{X},{\rm et}*}
((\ref{eqn:aclkumzn}))}}>> 
R\eta_{\os{\circ}{X},{\rm et}*}
(R^k\eps_{{\rm an}*}(\eps_{\rm an}^{-1}
\eta^{-1}_{\os{\circ}{X},{\rm et}}(E))) \\
@A{}AA @AA{(\ref{eqn:agnrtf})}A \\
\bigwedge^k({\cal M}^{\rm gp}_{X}
/{\cal O}^*_{X})\otimes_{\mab Z}E(-k) 
@>{\us{\sim}{(\ref{cd:knis})}}>> 
R^k\eps_{{\rm et}*}(\eps^{-1}_{\rm et}(E)).  
\end{CD}
\tag{5.4.1}\label{cd:kuagcrent}
\end{equation*}   
Furthermore, if $E$ is constructible, then 
the left vertical morphism is an isomorphism.
\end{prop}
\begin{proof}
The proof of the commutativity of (\ref{cd:kuagcrent}) 
is the same as that of the commutativity of (\ref{cd:kugcret}).  
By (\ref{theo:up}) for $\os{\circ}{X}$, 
the left vertical morphism is an isomorphism 
if $E$ is constructible.
\end{proof}

\section{Log exponential sequences}\label{sec:les} 
Let the notations be as in \S\ref{sec:lgelk}. 
Let $Y$ be an fs log analytic space 
over ${\mab C}$. 
In this section, 
introducing a new abelian sheaf 
${\cal L}^{\dag}_{Y^{\log}}$ 
in $\wt{Y}^{\log}$ 
which is a variant of the sheaf of logarithms 
${\cal L}_{Y^{\log}}$ 
defined in \cite{kn} 
and using the commutative diagram (\ref{cd:kugcret}), 
we give a commutative diagram 
(\ref{cd:cogokugcret}) below which 
compares the calculations of certain higher direct images 
by the use of the log exponential sequence in [loc.~cit.] 
and by the use of 
the algebraic log Kummer sequence in [loc.~cit.].  
We also introduce 
a new log exponential sequence 
which corrects 
the commutative diagram in \cite[(5.9.1)]{illl} 
(see (\ref{eqn:zlmzmlm}) and (\ref{rema:dcd}) (1) below). 
\par 
First let us recall the sheaf 
${\cal L}_{Y^{\log}}$ of logarithms of local sections 
of $\eps^{-1}_{\rm cl}({\cal M}^{\rm gp}_Y)$ 
in \cite[(1.4)]{kn}. 
\par 
Let 
${\rm Cont}_{Y^{\log}}(~, T)$ 
be a sheaf in 
$\wt{Y}^{\log}$ of continuous functions 
to a commutative topological group $T$.  
For a morphism $S \lo T$ of commutative topological groups,    
we have a natural morphism  
${\rm Cont}_{Y^{\log}}(~, S) \lo 
{\rm Cont}_{Y^{\log}}(~, T)$ of abelian sheaves in 
$\wt{Y}^{\log}$.  
The sheaf ${\cal L}_{Y^{\log}}$ is, by definition, 
the following fiber product 
\begin{equation*} 
{\cal L}_{Y^{\log}}:= 
{\rm Cont}_{Y^{\log}}(~, \sqrt{-1}{\mab R})
\times_{\exp,{\rm Cont}_{Y^{\log}}(~, {\mab S}^1)}
\eps^{-1}_{\rm cl}({\cal M}^{\rm gp}_{Y}).  
\tag{6.0.1}\label{eqn:ldefm}
\end{equation*}  
Then we have an exponential sequence 
\begin{equation*}
0 \lo {\mab Z}(1) \lo {\cal L}_{Y^{\log}} 
\os{\exp}{\lo} 
\eps^{-1}_{\rm cl}({\cal M}^{\rm gp}_{Y}) 
\lo 0  
\tag{6.0.2}\label{eqn:kmlmcl}
\end{equation*} 
in $\wt{Y}^{\log}$ ([loc.~cit., (1.4)]). 
\par 
By \cite[(1.5)]{kn}, for an abelian sheaf $E$ on 
$\wt{\os{\circ}{Y}}$, we have a canonical isomorphism 
\begin{equation*}
\bigwedge^k({\cal M}^{\rm gp}_{Y}/
{\cal O}^*_{Y})(-k)\otimes_{\mab Z}E
\os{\sim}{\lo}
R^k\eps_{{\rm cl}*}(\eps^{-1}_{\rm cl}(E)) \quad 
(k \in {\mab Z}_{\geq 0}).   
\tag{6.0.3}\label{eqn:mzcles}
\end{equation*}  
\par 
Next let us define an abelian sheaf 
${\cal L}^{\dag}_{Y^{\log}}$ 
in $\wt{Y}^{\log}$.  
\par 

We have a natural morphism 
\begin{equation*} 
\be^{-1}({\cal M}^{\rm gp}_{Y,\log}) 
\lo {\rm Cont}_{Y^{\log}}
(~, {\mab S}^1) 
\tag{6.0.4}\label{eqn:mcmmda}
\end{equation*} 
of abelian sheaves in $\wt{Y}^{\log}$ 
induced by the natural morphism 
\begin{equation*}
\Gam(V,{\cal M}^{\rm gp}_{Y,\log}) \lo 
{\rm Cont}_{Y^{\log}}
(V^{\log}, {\mab S}^1) \quad (V \in Y^{\log}_{\rm et})
\tag{6.0.5}\label{eqn:gvmcvs}
\end{equation*}
of presheaves on $Y^{\log}$.

\parno
Set
\begin{equation*} 
{\cal L}^{\dag}_{Y^{\log}}:= 
{\rm Cont}_{Y^{\log}}
(~, \sqrt{-1}{\mab R})
\times_{\exp,{\rm Cont}_{Y^{\log}}
(~, {\mab S}^1)}
\be^{-1}({\cal M}^{\rm gp}_{Y,\log}). 
\tag{6.0.6}\label{eqn:aldefm}
\end{equation*} 
Then we have an exponential sequence 
\begin{equation*}
0 \lo {\mab Z}(1) \lo {\cal L}^{\dag}_{Y^{\log}}
\os{\exp}{\lo} 
\be^{-1}({\cal M}^{\rm gp}_{Y,\log}) 
\lo 0.  
\tag{6.0.7}\label{eqn:akmlmcl}
\end{equation*} 

\begin{defi}
We call ${\cal L}^{\dag}_{Y^{\log}}$ 
the {\it sheaf of logarithms of local sections} of 
$\be^{-1}({\cal M}^{\rm gp}_{Y,\log})$. 
\end{defi}

\begin{lemm}\label{lemm:cbatsu}
The natural composite morphism 
\begin{equation*}
{\mab C}^* \os{\sus}{\lo} 
\be^{-1}({\cal M}^{\rm gp}_{Y,\log}) 
\lo {\rm Cont}_{Y^{\log}}
(~, {\mab S}^1) 
\tag{6.2.1}\label{eqn:umzczs}
\end{equation*} 
of abelian sheaves in $\wt{Y}^{\log}$ 
is induced by the map  
$c\lom c/\vert c \vert$ $(c \in {\mab C}^*)$.  
\end{lemm}
\begin{proof} 
The proof of (\ref{lemm:cbatsu}) is clear because, for any 
object $V$ of $Y^{\log}_{\rm et}$, 
${\mab C}^* \subset {\cal O}^*_{V,x}$ 
for any point $x$ of $\os{\circ}{V}$.  
\end{proof}

\begin{lemm}\label{lemm:mlzcl} 
For a positive integer $m$,  
the multiplication morphism 
\begin{equation*}
m\times \col 
{\cal L}^{\dag}_{Y^{\log}} 
\lo 
{\cal L}^{\dag}_{Y^{\log}} 
\tag{6.3.1}\label{eqn:msurj}
\end{equation*} 
is an isomorphism. 
\end{lemm}
\begin{proof} 
First we show the injectivity of (\ref{eqn:msurj}). 
Let $(a,s)$ 
$(a \in \sqrt{-1}{\mab R}, 
s\in \be^{-1}({\cal M}^{\rm gp}_{Y,\log}))$ 
be a local section of 
${\cal L}^{\dag}_{Y^{\log}}$ 
such that $m(a,s)=0$. 
Then $a=0$ and $s\in ({\mab Z}/m)(1)$. 
By (\ref{lemm:cbatsu}), we see that $s=1$.  
Hence the morphism (\ref{eqn:msurj}) is injective. 
\par 
Next we show the surjectivity of (\ref{eqn:msurj}). 
For an object $U$ of $Y^{\log}$, 
let $(a,u)$ 
$(a \in \sqrt{-1}{\mab R}, 
u\in \Gam(U,\be^{-1}({\cal M}^{\rm gp}_{Y,\log})))$ 
be a section of 
$\Gam(U,{\cal L}^{\dag}_{Y^{\log}})$. 
We may assume that 
$u \in \Gam(U,\be^{-1}({\cal M}_{Y,\log}))$. 
Since ${\cal M}_{Y,\log}$ is $m$-divisible 
by (\ref{lemm:logkum}) 
and since the functor $\be^{-1}$ is right-exact, 
there exists a section $v_{\lam}$ of 
$\Gam(U_{\lam}, \be^{-1}({\cal M}_{Y,\log}))$ for  
some covering $(U_{\lam} \lo U)_{\lam}$ of $U$ 
in $Y^{\log}$ 
such that $v^m_{\lam}=u\vert_{U_{\lam}}$. 
Let $w_{\lam}$ be the image of $v_{\lam}$ 
in 
${\rm Cont}_{Y^{\log}}(U_{\lam}, {\mab S}^1)$. 
Then $\zeta_{\lam}:=w_{\lam}\exp(-m^{-1}a)$ is 
an $m$-th root of unity. 
Hence $(m^{-1}a, v_{\lam}\zeta^{-1}_{\lam})$ is 
indeed an element of 
$\Gam(U_{\lam},{\cal L}^{\dag}_{Y^{\log}})$ by 
(\ref{lemm:cbatsu}), and  
$m(m^{-1}a,v_{\lam}\zeta^{-1}_{\lam})=(a,u)\vert_{U_{\lam}}$. 
Hence the morphism (\ref{eqn:msurj}) is surjective. 
\end{proof}

\begin{prop}\label{coro:fzcl}
Let $E$ be an abelian sheaf in $\wt{\os{\circ}{Y}}$. 
Then the following diagram 
\begin{equation*}
\begin{CD}  
\bigwedge^k
({\cal M}^{\rm gp}_{Y}/{\cal O}_{Y}^*)(-k)
\otimes_{\mab Z}E 
@>{\us{\sim}{(\ref{eqn:mzcles})}}>> 
R^k\eps_{{\rm cl}*}(\eps^{-1}_{\rm cl}(E))  \\ 
@V{{\rm id}\otimes{\exp(m^{-1}\times)}^{\otimes k}
\otimes{\rm proj}}VV 
@VVV \\ 
\bigwedge^k ({\cal M}^{\rm gp}_{Y}/
{\cal O}_{Y}^*)\otimes_{\mab Z}({\mab Z}/m)(-k)
\otimes_{\mab Z}(E/mE)   
@>{\us{\sim}{(\ref{eqn:clkumz})}}>>  
R^k\eps_{{\rm cl}*}(\eps^{-1}_{\rm cl}(E/mE)) 
\end{CD} 
\tag{6.4.1}\label{eqn:anmoclm}
\end{equation*} 
is commutative for $k\in {\mab Z}_{\geq 0}$.  
\end{prop}
\begin{proof}
By (\ref{lemm:mlzcl}) 
we obtain a well-defined morphism 
\begin{equation*}
\exp(m^{-1}\times)  \col 
{\cal L}^{\dag}_{Y^{\log}} \lo 
\be^{-1}({\cal M}^{\rm gp}_{Y,\log}).  
\tag{6.4.2}\label{eqn:expmm}
\end{equation*} 
Since $(2\pi \sqrt{-1}n/m, \exp(2\pi \sqrt{-1}n/m))$ 
$(n\in {\mab Z})$ 
is a section of 
${\cal L}^{\dag}_{Y^{\log}}$ 
by (\ref{lemm:cbatsu}), 
$$m^{-1}(2\pi \sqrt{-1}n, 1)=
(2\pi \sqrt{-1}n/m, \exp(2\pi \sqrt{-1}n/m)).$$ 
Hence we obtain the following commutative diagram 
\begin{equation*}
\begin{CD} 
0 @>>> {\mab Z}(1) @>>> 
{\cal L}_{Y^{\log}} 
@>{\exp}>> 
\eps^{-1}_{\rm cl}({\cal M}^{\rm gp}_{Y})
@>>> 0 \\
@. @| @VVV @VV{(\ref{eqn:zmbe})^{\rm gp}}V \\ 
0 @>>> {\mab Z}(1) @>>> 
{\cal L}^{\dag}_{Y^{\log}} 
@>{\exp}>> 
\be^{-1}({\cal M}^{\rm gp}_{Y,\log}) 
@>>> 0  \\
@. @V{\exp(m^{-1}\times)}VV 
@V{\exp(m^{-1}\times)}VV @| @. \\
0 @>>> ({\mab Z}/m)(1) @>>> 
\be^{-1}({\cal M}^{\rm gp}_{Y,\log}) 
@>{m\times}>> 
\be^{-1}({\cal M}^{\rm gp}_{Y,\log}) @>>> 0   
\end{CD}
\tag{6.4.3}\label{cd:lmuek}
\end{equation*} 
of exact sequences.  
(\ref{coro:fzcl}) immediately follows from 
the commutative diagram (\ref{cd:lmuek}) and 
from the definitions of the isomorphisms 
(\ref{eqn:mzcles}) and  (\ref{eqn:clkumz}). 
\end{proof}

\begin{coro}\label{coro:final} 
Let $X$ be an fs log scheme over ${\mab C}$ 
whose underlying scheme $\os{\circ}{X}$ is 
locally of finite type over ${\mab C}$. 
Let $m$ be a positive integer. 
Let $E$ be an $m$-torsion abelian sheaf 
in $\wt{\os{\circ}{X}}_{\rm et}$. 
Then the following diagram is commutative$:$
\begin{equation*} 
\begin{CD}
R\eta_*(\bigwedge^k(\eta^*({\cal M}^{\rm gp}_X)
/{\cal O}^*_{X_{\rm an}})\otimes_{\mab Z}\eta^{-1}(E)(-k))
@>{\us{\sim}{R\eta_*((\ref{eqn:mzcles}))}}>> 
R\eta_*
(R^k\eps_{{\rm cl}*}(\eps_{\rm cl}^{-1}\eta^{-1}(E))) \\ 
@V{R\eta_*({\rm id}\otimes 
{\rm id}\otimes\exp(m^{-1}\times))}V{\simeq}V 
@| \\ 
R\eta_*(\bigwedge^k(\eta^*({\cal M}^{\rm gp}_X) 
/{\cal O}^*_{X_{\rm an}})\otimes_{\mab Z}\eta^{-1}(E)(-k))
@>{\us{\sim}{R\eta_*((\ref{eqn:clkumz}))}}>> 
R\eta_*(R^k\eps_{{\rm cl}*}
(\eps_{\rm cl}^{-1}\eta^{-1}(E))) \\
@A{}AA @AAA \\
\bigwedge^k({\cal M}^{\rm gp}_{X}
/{\cal O}^*_{X})\otimes_{\mab Z}E(-k) 
@>{\us{\sim}{(\ref{cd:knis})}}>> 
R^k\eps_{{\rm et}*}(\eps^{-1}_{\rm et}(E)).  
\end{CD}
\tag{6.5.1}\label{cd:cogokugcret}
\end{equation*}  
\end{coro} 
\begin{proof} 
By (\ref{coro:fzcl}) 
we obtain the upper commutative diagram; the lower one 
is nothing but the commutative diagram (\ref{cd:kugcret}).
\end{proof}


Lastly we point out the mistakes in 
the proof of \cite[(5.9)]{illl} 
and we correct them. 
\par 
In the proof of \cite[(5.9)]{illl}, it is claimed 
that ${\cal L}_{X^{\log}_{\rm an}}$  is 
uniquely divisible by 
a positive integer $m$ for an fs log scheme $X$ 
over ${\mab C}$ whose underlying scheme is 
locally of finite type over ${\mab C}$. 
However the divisibility 
does not hold for $m \geq 2$. 
Indeed, if  ${\cal L}_{X^{\log}_{\rm an}}$ were  
$m$-divisible, then 
$\eps^{-1}_{\rm cl}({\cal M}^{\rm gp}_{X_{\rm an}})$ 
would also be by the exponential sequence 
(\ref{eqn:kmlmcl}). 
Consequently 
$\eps^{-1}_{\rm cl}
({\cal M}^{\rm gp}_{X_{\rm an}}/{\cal O}^*_{X_{\rm an}})$
would also be $m$-divisible. 
Clearly this does not hold in general. 
Indeed, 
let $X$ be the log point $s=({\rm Spec}\,{\mab C}, 
{\mab N}\oplus{\mab C}^*)$. 
Then 
$\eps^{-1}_{\rm cl}
({\cal M}^{\rm gp}_{X_{\rm an}}/
{\cal O}^*_{X_{\rm an}})=
\eps^{-1}_{\rm cl}({\mab Z})={\mab Z}$. 
Therefore the proof of \cite[(5.9)]{illl} 
is mistaken. 
\par 
I think that there does not exist 
the natural morphism 
$\tau^{-1}(M^{\rm gp}) 
\lo \eta^{-1}(M^{\rm gp}_{X_{\rm ket}})$ 
in the commutative diagram \cite[(5.9.1)]{illl} 
(I think that there does not exist 
a useful direct relation between 
the log exponential sequence 
in \cite{kn} and the algebraic log Kummer sequence in 
[loc.~cit.]): 
consider the trivial log case. Then $M^{\rm gp}$ 
in \cite[(5.9.1)]{illl} 
is the multiplicative group 
of the sheaf of germs of invertible 
holomorphic functions on $X_{\rm an}$ 
(not the sheaf of invertible 
algebraic functions on $X$) 
and the morphism $\tau^{-1}(M^{\rm gp}) 
\lo \eta^{-1}(M^{\rm gp}_{X_{\rm ket}})$ 
is a morphism from an analytic sheaf 
to an algebraic sheaf.  
Usually such a morphism does not naturally 
exist except the trivial morphism.

\par
Let us give a right commutative diagram in 
(\ref{eqn:zlmzmlm}) below 
(see also (\ref{rema:dcd}) (1) below). 
\par
Let the notations  be as in \S\ref{sec:knct}. 
Then we have a natural morphism 
\begin{equation*} 
\eta^{{\log},-1}({\cal M}^{\rm gp}_{X,\log}) 
\lo {\rm Cont}_{X^{\log}_{\rm an}}
(~, {\mab S}^1) 
\tag{6.5.2}\label{eqn:mcmm}
\end{equation*} 
induced by the natural morphisms 
$\Gam(U,{\cal M}^{\rm gp}_{X,\log}) \lo 
{\rm Cont}_{X^{\log}_{\rm an}}
(U^{\log}_{\rm an}, {\mab S}^1)$ 
of abelian presheaves for objects $U$'s of 
$X^{\log}_{\rm et}$. 
Set
\begin{equation*} 
{\cal L}^{-1}_{X^{\log}_{\rm an}}= 
{\rm Cont}_{X^{\log}_{\rm an}}
(~, \sqrt{-1}{\mab R})
\times_{\exp,{\rm Cont}_{X^{\log}_{\rm an}}
(~, {\mab S}^1)}
\eta^{\log,-1}({\cal M}^{\rm gp}_{X,\log}). 
\tag{6.5.3}\label{eqn:ldefmna}
\end{equation*} 
Then we have an exponential sequence 
\begin{equation*}
0 \lo {\mab Z}(1) \lo 
{\cal L}^{-1}_{X^{\log}_{\rm an}}
\os{\exp}{\lo} 
\eta^{\log,-1}({\cal M}^{\rm gp}_{X,\log}) 
\lo 0.  
\tag{6.5.4}\label{eqn:kmlmclna}
\end{equation*}

\par
Let us also recall 
$\eta^{\log *}({\cal M}^{\rm gp}_{X,\log})$ in 
(\ref{eqn:stea}).  
We have a natural morphism 
\begin{equation*} 
\eta^{\log *}({\cal M}^{\rm gp}_{X,\log}) 
\lo {\rm Cont}_{X^{\log}_{\rm an}}
(~, {\mab S}^1) 
\tag{6.5.5}\label{eqn:mcmmst}
\end{equation*} 
induced by the natural morphism  
$$\Gam(U,{\cal M}^{\rm gp}_{X,\log}) \lo 
{\rm Cont}_{X^{\log}_{\rm an}}
(U^{\log}_{\rm an}, {\mab S}^1)$$ 
of abelian presheaves for objects $U$'s of 
$X^{\log}_{\rm et}$ 
and by the following morphism  
$$\Gam(V,{\cal O}^*_{X_{\rm an},\log})\owns f \lom 
((x,h)\lom f(x)/\vert f(x) \vert) \in 
{\rm Cont}_{X^{\log}_{\rm an}}
(V^{\log}, {\mab S}^1)$$ 
of abelian presheaves for objects $V$'s of 
$(X_{\rm an})^{\log}_{\rm et}$. 
Here $x$ is a point of $\os{\circ}{V}$ 
and $h \col {\cal M}^{\rm gp}_{V,x} \lo {\mab S}^1$ 
is a morphism of groups such 
that $h(g)=g(x)/\vert g(x) \vert$ for any 
$g \in {\cal O}^*_{V,x}$. 
Set 
\begin{equation*} 
{\cal L}^{*}_{X^{\log}_{\rm an}}:= 
{\rm Cont}_{X^{\log}_{\rm an}}
(~, \sqrt{-1}{\mab R})
\times_{\exp,{\rm Cont}_{X^{\log}_{\rm an}}
(~, {\mab S}^1)}
\eta^{\log *}({\cal M}^{\rm gp}_{X,\log}). 
\tag{6.5.6}\label{eqn:aldefmst}
\end{equation*} 
Then we have an exponential sequence 
\begin{equation*}
0 \lo {\mab Z}(1) \lo 
{\cal L}^{*}_{X^{\log}_{\rm an}}
\os{\exp}{\lo} 
\eta^{\log *}({\cal M}^{\rm gp}_{X,\log}) 
\lo 0.  
\tag{6.5.7}\label{eqn:akmlmstl}
\end{equation*}

\begin{defi}
Let $\natural$ be $-1$ or $*$. 
We call ${\cal L}^{\natural}_{X^{\log}_{\rm an}}$
the {\it sheaf of logarithms of local sections} of 
$\eta^{\log \natural}({\cal M}^{\rm gp}_{X,\log})$. 
\end{defi}

\begin{lemm}\label{lemm:ldiv-1} 
$(1)$ The natural composite morphism 
\begin{equation*}
{\mab C}^* \os{\sus}{\lo} 
\eta^{\log \natural}({\cal M}^{\rm gp}_{X,\log}) 
\lo {\rm Cont}_{X^{\log}_{\rm an}}
(~, {\mab S}^1) 
\tag{6.7.1}\label{eqn:-1umzs}
\end{equation*} 
of abelian sheaves in $\wt{X}^{\log}_{\rm an}$ 
is induced by the map  
$c\lom c/\vert c \vert$ $(c \in {\mab C}^*)$.  
\par 
$(2)$ For a positive integer $m$, 
the multiplication morphism 
\begin{equation*}
m\times \col 
{\cal L}^{\natural}_{X^{\log}_{\rm an}} 
\lo 
{\cal L}^{\natural}_{X^{\log}_{\rm an}} 
\tag{6.7.2}\label{eqn:msurjcl}
\end{equation*} 
is an isomorphism. 
\end{lemm}
\begin{proof} 
The proof of (1) (resp.~(2)) is the same as that of 
(\ref{lemm:cbatsu}) (resp.~(\ref{lemm:mlzcl})). 
\end{proof} 

By (\ref{lemm:ldiv-1}) (2)
we obtain a well-defined morphism 
\begin{equation*}
\exp(m^{-1}\times)  \col 
{\cal L}^{\natural}_{X^{\log}_{\rm an}} \lo 
\eta^{\log \natural}({\cal M}^{\rm gp}_{X,\log})   
\tag{6.7.3}\label{eqn:expmml}
\end{equation*} 
and the following commutative diagram 
\begin{equation*}
\begin{CD} 
0 @>>> {\mab Z}(1) @>>> 
{\cal L}^{\natural}_{X^{\log}_{\rm an}} 
@>{\exp}>> 
\eta^{\log \natural}({\cal M}^{\rm gp}_{X,\log}) 
@>>> 0  \\
@. @V{\exp(m^{-1}\times)}VV 
@V{\exp(m^{-1}\times)}VV @| @. \\
0 @>>> ({\mab Z}/m)(1) @>>> 
\eta^{\log \natural}({\cal M}^{\rm gp}_{X,\log})  
@>{m\times}>> \eta^{\log \natural}({\cal M}^{\rm gp}_{X,\log}) 
@>>> 0   
\end{CD}
\tag{6.7.4}\label{cd:nlmuek}
\end{equation*} 
of exact sequences. 

As a summary,  we obtain the following: 
\begin{prop}\label{prop:bgc} 
There exists the following commutative diagram 
\begin{equation*}
\begin{CD} 
0 @>>> {\mab Z}(1) @>>> 
{\cal L}_{X^{\log}_{\rm an}} 
@>{\exp}>> 
\eps^{-1}_{\rm cl}(\eta^*({\cal M}^{\rm gp}_X)) 
@>>> 0 \\
@. @| @VVV @VV{(\ref{eqn:ster})^{\rm gp}}V \\ 
0 @>>> {\mab Z}(1) @>>> 
{\cal L}^{\dag}_{X^{\log}_{\rm an}} 
@>{\exp}>> 
\be^{-1}({\cal M}^{\rm gp}_{X_{\rm an},\log}) 
@>>> 0 \\
@. @| @AAA @AA{\be^{-1}((\ref{eqn:etmx})^{\rm gp})}A \\
0 @>>> {\mab Z}(1) @>>> 
{\cal L}^{*}_{X^{\log}_{\rm an}} 
@>{\exp}>> 
\eta^{\log *}({\cal M}^{\rm gp}_{X,\log}) 
@>>> 0 \\ 
@. @| @AAA @AA{(\ref{eqn:-1*et})^{\rm gp}}A \\
0 @>>> {\mab Z}(1) @>>> 
{\cal L}^{-1}_{X^{\log}_{\rm an}} 
@>{\exp}>> 
\eta^{\log,-1}({\cal M}^{\rm gp}_{X,\log}) 
@>>> 0 \\ 
@. @V{\exp(m^{-1}\times)}VV 
@V{\exp(m^{-1}\times)}VV  @| @. \\
0 @>>> ({\mab Z}/m)(1) @>>> 
\eta^{\log,-1}({\cal M}^{\rm gp}_{X,\log})  
@>{m\times}>>  
\eta^{\log,-1}({\cal M}^{\rm gp}_{X,\log}) @>>> 0   
\end{CD}
\tag{6.8.1}\label{eqn:zlmzmlm}
\end{equation*}
of exact sequences.
\end{prop}  

\begin{rema}\label{rema:dcd} 
(1) Using the morphism (\ref{eqn:stepet}), we can delete 
the second exact sequence in (\ref{eqn:zlmzmlm}). 
\par
(2) Using the commutative diagram 
(\ref{eqn:zlmzmlm}) (or the commutative diagram obtained in (1)), 
we can give a proof of \cite[(2.6)]{kn}. 
\end{rema}

\parno
Department of Mathematics, Tokyo Denki University,
2--2 Kanda Nishiki-cho Chiyoda-ku, Tokyo 101--8457, Japan. 
\parno
{\it E-mail address\/}: nakayuki@cck.dendai.ac.jp
\end{document}